\documentclass{elsarticle}

\usepackage{hyperref}
\usepackage{amsmath,amssymb}
\usepackage{tikz}

\hypersetup{pdfborder=0 0 0}

\newcommand{\V}{V}

\newcommand{\Vn}[1]{V_{#1}}
\newcommand{\Vplus}{V_+}
\newcommand{\VR}[1]{V^{#1}}
\newcommand{\un}[1]{u_{#1}}
\newcommand{\uR}[1]{u^{#1}}

\newcommand{\Un}[1]{\mathcal U_{#1}}
\newcommand{\Umn}[2]{\mathcal U^{#1}_{#2}}

\newcommand{\var}[1]{\vec a^{#1}}
\newcommand{\an}[1]{a_{#1}}
\newcommand{\arn}[2]{a^{#1}_{#2}}
\newcommand{\br}[1]{b^{#1}}
\newcommand{\massMatrix}{\mathcal M}
\newcommand{\snapshotMatrix}{\mathcal U}
\newcommand{\snapshotGramian}{\mathcal G}
\newcommand{\Rmax}{D}
\newcommand{\parameter}{\mu}
\newcommand{\parameterDomain}{S}
\newcommand{\iParam}{p}
\newcommand{\nParam}{P}
\newcommand{\iTime}{k}
\newcommand{\nTime}{K}
\newcommand{\iSnap}{n}
\newcommand{\nSnap}{N}

\newcommand{\nPod}{R}

\DeclareMathOperator{\spn}{span}

\begin{document}

\begin{frontmatter}

\title{POD-Galerkin reduced-order modeling with \\adaptive finite element snapshots}

\author[mymainaddress,mysecondaryaddress]{Sebastian Ullmann}
\author[mysecondaryaddress]{Marko Rotkvic}
\author[mymainaddress,mysecondaryaddress,mythirdaddress]{Jens Lang}
\address[mymainaddress]{Graduate School Computational Engineering, Technische Universit\"at Darmstadt, \\Dolivostr.\ 15, D-64293 Darmstadt, Germany}
\address[mysecondaryaddress]{Department of Mathematics, Technische Universit\"at Darmstadt, \\Dolivostr.\ 15, D-64293 Darmstadt, Germany}
\address[mythirdaddress]{Graduate School Energy Science and Engineering, Technische Universit\"at Darmstadt,\\Jovanka-Bontschits-Str.\ 2, 64287 Darmstadt}

\begin{abstract}
We consider model order reduction by proper orthogonal decomposition (POD) for parametrized partial differential equations, where the underlying snapshots are computed with adaptive finite elements. We address computational and theoretical issues arising from the fact that the snapshots are members of different finite element spaces. We propose a method to create a POD-Galerkin model without interpolating the snapshots onto their common finite element mesh. The error of the reduced-order solution is not necessarily Galerkin orthogonal to the reduced space created from space-adapted snapshot. We analyze how this influences the error assessment for POD-Galerkin models of linear elliptic boundary value problems. As a numerical example we consider a two-dimensional convection-diffusion equation with a parametrized convective direction. To illustrate the applicability of our techniques to non-linear time-dependent problems, we present a test case of a two-dimensional viscous Burgers equation with parametrized initial data.
\end{abstract}


\end{frontmatter}

\section{Introduction}

Model order reduction is a tool to decrease the computational cost for applications where a pa\-ra\-me\-tri\-zed PDE problem needs to be solved multiple times for different parameter values. Therefore model order reduction is often studied in the context of optimal control \cite{GhiglieriUlbrich2014,NegriEA2013,TroltzschVolkwein2009} or uncertainty quantification \cite{ChenQuarteroni2015,HaasdonkEA2013,UllmannLang2014}. Snapshot-based model order reduction requires a set of representative samples of the solution, which need to be computed in advance. The solution of the reduced-order model is then represented as a linear combination of these snapshots. The respective coefficients are determined by means of a Galerkin projection, based on a weak form of the governing equations. In this way, the reduced-order model inherits both the spatial structure of typical solutions as well as the underlying physics. Introductions to snapshot-based model order reduction are provided by the textbooks \cite{HesthavenEA2016,QuarteroniEA2016}.

Standard techniques for model order reduction assume that all snapshots use one and the same spatial mesh. We refer to this case as \emph{static} snapshot computations. In contrast, with \emph{adaptive} snapshot computations we mean that each snapshot may use a different mesh. Combining space-adaptive simulations with model order reduction can be an advantage from two points of view: Firstly, introducing spatial adaptivity to the set-up phase of a reduced-order model can decrease the total computation time if the solution contains local features depending on the parameter. By adapting the mesh to the features at a given parameter value, less degrees of freedom are required to obtain a certain accuracy. Secondly, introducing model order reduction to space-adaptive simulations promises computational speed-up in cases where the solution is susceptible to an approximation in a low-dimensional linear space and needs to be evaluated for multiple different parameter values.

One possible way to implement a POD-Galerkin reduced-order model from snapshots of different discretization spaces is to express the snapshots as elements of some common discretization space. Then one can use standard methods to create a POD basis and compute a respective Galerkin projection of the solution. In the case of strongly varying local refinements, however, a good common discretization space may be relatively high-dimensional, which makes it unattractive from a computational point of view. We show that by expressing the POD basis in terms of the snapshots, it is possible to avoid forming the common discretization space explicitly.

Model order reduction with spatial adaptivity has been studied in \cite{AliUrbanSteih20XX,AliUrban2016} for snapshot computations with adaptive wavelets and in \cite{Yano2016} for snapshot computations with adaptive mixed finite elements. The main issue addressed in these publications is the assessment of the error between the reduced-order solution and the infinite-dimensional true solution. In the case of static snapshot computations, this problem can be circumvented by assuming a sufficiently fine snapshot discretization. Then the error between the reduced-order solution and a corresponding discrete solution can be estimated with the help of the discrete residual. For the case of adaptive snapshot computations, \cite{AliUrbanSteih20XX,AliUrban2016} use wavelet techniques to estimate the required dual norm of the continuous residual. In contrast, \cite{Yano2016} derives a bound for the dual norm of the continuous residual from a special mixed finite element and reduced basis formulation.

The references \cite{AliUrbanSteih20XX,AliUrban2016,Yano2016} focus on model order reduction by the greedy reduced basis method \cite{PrudhommeEA2002}. In this paper, however, we consider an alternative approach, namely proper orthogonal decomposition (POD) \cite{HolmesEA1996,Sirovich1987}. The major difference between both methods lies in the construction of the reduced space used as a test and trial space in a Galerkin procedure. Both methods require a fixed set of training parameters to be chosen in advance, where for time dependent problems, time is viewed as a parameter. From the training parameter set, the greedy reduced basis method selects a set of parameter values in an iterative way using an error estimator and uses the span of the corresponding snapshots as a reduced space. This can save computation time by avoiding the computation of snapshots which have not been selected. The resulting reduced space is close to the one which minimizes the maximum approximation error over the training set. In contrast, POD forms a reduced basis by linearly combining the snapshots corresponding to all training parameter values. The linear combination is done in a way which minimizes the mean square approximation error over the training set, but at the cost of computing all training snapshots. If the dimension of the reduced spaces are increased, both the greedy and the POD space eventually become equal to the span of the snapshots corresponding to all training parameter values. 

A POD of snapshots resulting from a static spatial discretization can be implemented in terms of a truncated singular value decomposition \cite{KunischVolkwein1999}. Therefore, an error estimator is not necessary for creating a POD reduced basis. Because snapshots have to be computed for all training parameters, POD is often applied to time-dependent problems, where snapshot data arise as a by-product of the numerical time stepping scheme. We note that in this context, POD with time-adaptive snapshots has been studied in \cite{AllaEA2015}. POD with one-dimensional space-adaptive snapshots has already been addressed in \cite{Lass2014}, where the POD computation relies on a polynomial approximation of the snapshots. In contrast, we focus on the two-dimensional case and present a method which does not require an intermediate approximation of the snapshots.

A major difference between greedy and POD reduced basis methods for space adaptive snapshots is caused by the relation between the snapshots and the reduced basis functions: In the greedy reduced basis method, the reduced space is formed by linear combinations of snapshots, while the snapshots are themselves elements of the reduced space. In the POD reduced basis method, the reduced space is also formed by linear combinations of snapshots, but the snapshots are not elements of the reduced space, in general. The difference between the snapshots and their closest approximation in the POD space can be measured in terms of the truncated singular values. This has consequences for the error assessment of POD-Galerkin schemes in presence of space-adapted snapshots. While the main difficulties in the greedy reduced basis method arises from the fact that the error of the reduced-order solution is not necessarily orthogonal to the reduced space anymore, the POD reduced basis method is additionally subject to a truncation error.

This paper is structured as follows: In section \ref{sec:Pod}, we introduce proper orthogonal decomposition for adaptive finite element snapshots. We propose methods to efficiently compute POD bases for adaptive finite element discretizations with nested refinement. A POD-Galerkin reduced-order model based on adaptive snapshots is formulated in section \ref{sec:elliptic} for an elliptic boundary value problem. We prove error statements for the reduced-order solution in presence of adaptive snapshots and compare the results to snapshot computations on a static mesh. The methods and analytic results are illustrated in section \ref{sec:ConvectionDiffusion} with a numerical test case involving a linear convection-diffusion equation with parametrized convective direction. The applicability to non-linear time-dependent problems is suggested by the results of section \ref{sec:Burgers}, which features a Burgers problem with parametrized initial condition.

\section{Proper orthogonal decomposition}\label{sec:Pod}

We consider snapshot-based model order reduction, where the solution of a PDE problem is represented in the space spanned by a set of reduced basis functions obtained by linearly combining a set of snapshots. Such reduced-basis functions typically have a global support and contain information about expected spatial structures of the solution. 

One method to compute reduced basis functions from snapshots is the proper orthogonal decomposition \cite{HolmesEA1996,Sirovich1987}. If the snapshots correspond to coefficient vectors of a finite element discretization on a fixed grid, a POD can be given in terms of the snapshot coefficient matrix \cite{KunischVolkwein1999}. In the following, we introduce a method that computes a POD of snapshots stemming from an adaptive finite element simulation, where a snapshot matrix can not be created in straight-forward manner.

\subsection{Method of snapshots}\label{sec:MethodOfSnapshots}

We consider a PDE problem defined over some bounded open spatial domain $\Omega$ and some time and/or parameter domain $\parameterDomain$. In particular, we are interested in parametrized elliptic boundary value problems, where $\parameterDomain$ is a parameter domain, and in parametrized parabolic initial boundary value problems, where $\parameterDomain$ is a tensor product of a time interval and a parameter domain.

Let $\V$ be the infinite-dimensional Hilbert space used to characterize the solution as a function of space. A typical example is $V=H^1_0(\Omega)$, the Sobolev space of $L^2(\Omega)$ functions with weak first derivatives in $L^2(\Omega)$ and boundary values vanishing in the sense of traces. We denote the $V$-scalar product by $(\cdot,\cdot)_{\V}$ and the $V$-norm by $\|\cdot\|_{\V}$.

A proper orthogonal decomposition of snapshots $u_1,\dots,u_N\in V$ can be defined in terms of a system of minimization problems \cite{KunischVolkwein2001}: Find functions $\phi_1,\dots,\phi_N\in V$ which solve the minimization problems
\begin{align}\label{eq:PODminimization}
  &\min_{\phi_1,\dots,\phi_R\in V}\sum_{n=1}^N\Big\|u_n-\sum_{k=1}^R(u_n,\phi_k)_V\phi_k\Big\|^2_V,\quad(\phi_i,\phi_j)_V=\delta_{ij},\quad i,j=1,\dots,R
\end{align}
for all $R=1,\dots,N$. The solutions can be computed by an eigenvalue decomposition of the matrix containing the mutual $V$-inner products of $u_1,\dots,u_N$. This approach is often called the \emph{method of snapshots} \cite{Sirovich1987}. The eigenvalue problem can be written in components as follows:
For given $\un{1},\dots,\un{N}\in\V$, find $\lambda\in\mathbb R$ and $\vec a=(\an{1},\dots,\an{N})^T\in\mathbb R^{N}$, such that
\begin{align*}
  \sum_{j=1}^N(u_i,u_j)_{\V}\an{j}&=\lambda \an{i},\qquad i=1,\dots,N.
\end{align*}
After defining the snapshot Gramian matrix $\snapshotGramian = (g_{ij})$ with $g_{ij} = (u_i,u_j)_{\V}$ for $i,j=1,\dots,N$, the matrix form of the set of equations is given by $\snapshotGramian \vec a=\lambda \vec a$. The eigenvalue decomposition of the snapshot Gramian results in eigenvalues $\lambda_1,\dots,\lambda_{N}\in\mathbb R$ and eigenvectors $\var{1},\dots,\var{N}\in\mathbb R^N$. We order the eigenvalues such that $\lambda_1\geq\dots\geq\lambda_{\Rmax}>0 = \lambda_{\Rmax+1}=\dots=\lambda_{N}$ and write the eigenvectors in components as $\var{r}=(\arn{r}{1},\dots,\arn{r}{N})^T$ for $r=1,\dots,N$. Then the first POD basis functions are given by linear combinations of snapshots,
\begin{align}\label{eq:PODbasis}
  \phi_r=\sum_{n=1}^N \un{n}\frac{\arn{r}{n}}{\sqrt{\lambda_r}},\qquad r=1,\dots,\Rmax.
\end{align}
The space $\VR{R}=\spn(\phi_1,\dots,\phi_R)$ is called a POD space of dimension $R$ for any $1\leq R \leq D$.

POD reduced-order modeling tries to approximate a solution $u:\parameterDomain\rightarrow\V$ with a function $\uR{R}:\parameterDomain\rightarrow\VR{R}$ defined by
\begin{align}\label{eq:PODapproximation}
  \uR{R} = \sum_{r=1}^R\phi_r\br{r},
\end{align}
where $\vec b=(\br{1},\dots,\br{R})^T:\parameterDomain\rightarrow\mathbb R^R$ is a POD coefficient vector. Combining \eqref{eq:PODbasis} and \eqref{eq:PODapproximation} gives a POD approximation in terms of the snapshots,
\begin{align}\label{eq:PODapproximation2}
  \uR{R} = \sum_{n=1}^N \un{n}\sum_{r=1}^R\frac{\arn{r}{n}\br{r}}{\sqrt{\lambda_r}}.
\end{align}

One particular choice of POD coefficients is implied by the POD minimization problem \eqref{eq:PODminimization} and given by a $\V$-orthogonal projection of $u$ onto $\VR{R}$:
\begin{align*}
  P^Ru := \sum_{r=1}^R\phi_r(\phi_r,u)_V\quad\forall u\in\V,\quad R=1,\dots,D,
\end{align*}
which means $b^r=(\phi_r,u)_\V$ for $r=1,\dots,R$. The orthogonal projection is used as a reference solution later on, because it gives a POD representation with optimal coefficients:
\begin{align*}
  \|u-P^Ru\|_\V&=\inf_{v\in\VR{R}}\|u-v\|_\V\quad\forall u\in\V.
\end{align*}
The error of the POD projection of the snapshots can be computed from the POD eigenvalues,
\begin{align}\label{eq:sumEig}
  \sum_{n=1}^N \|u_n-P^Ru_n\|_V^2 = \sum_{n=R+1}^D \lambda_n,
\end{align}
which implies that the POD projection error of the snapshots decreases monotonically with the POD dimension and that $u_n=P^Du_n$ for $n=1,\dots,N$. Together with \eqref{eq:PODbasis} this means 
\begin{align}\label{eq:equalSpans}
  \spn(\un{1},\dots,\un{N})=\spn(\phi_1,\dots,\phi_{\Rmax}).
\end{align}

\subsection{Adaptive snapshot spaces}

In order to compute a set of snapshots, we discretize our PDE problem of interest with adaptive finite elements in space. Let $\Vn{1},\dots,\Vn{N}\subset\V$ be adapted finite element spaces, so that $\un{1}\in\Vn{1},\dots,\un{N}\in\Vn{N}$. Let $M_{1},\dots,M_{N}$ be the dimensions of the respective spaces. We focus on $h$-adaptive Lagrangian finite elements with a fixed polynomial degree, so that each snapshot finite element space is defined by a triangulation.

For discretizations on a fixed triangulation, one can represent linear combinations of snapshots by linear combinations of finite element coefficient vectors. In order to do this for adaptive spatial discretizations, however, one must first express the snapshots in terms of a suitable common finite element basis. Therefore, we introduce a space $\Vplus\subset V$ with finite dimension $M_+$, on which we impose two properties:
\begin{enumerate}
  \item $\Vplus$ is a finite element space of the same type as $\Vn{1},\dots,\Vn{N}$,
  \item $\Vn{1}+\dots+\Vn{N}\subset\Vplus$ in terms of a vector sum.
\end{enumerate}
A consequence of the first property is that after interpolating all snapshots onto $\Vplus$, we can work with them in the same way as if they were computed on a fixed triangulation. The second property ensures that the error between any snapshot and its representation in $\Vplus$ is zero. 

\begin{figure}
  \begin{center}
    \includegraphics{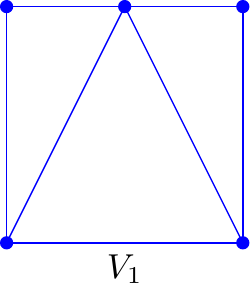}\quad
    \includegraphics{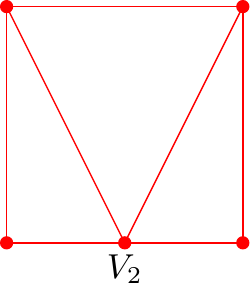}\quad
    \includegraphics{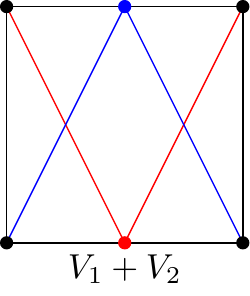}\quad
    \includegraphics{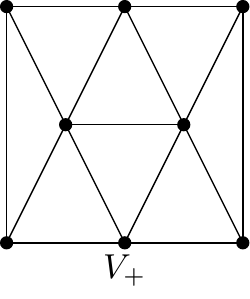}
  \end{center}
  \vspace{-5mm}
  \caption{Illustration of meshes corresponding to general finite element spaces $V_1$ and $V_2$, their vector sum $V_1+V_2$ and a common finite element space $V_+$ obtained by adding nodes and edges.}
  \label{fig:SketchFEUnion}
\end{figure}

In general, setting $\Vplus=\Vn{1}+\dots+\Vn{N}$ would be too restrictive in the sense that it does not necessarily fulfill the first property. Consider, for example, the case where $\Vn{1},\dots,\Vn{N}$ are linear Lagrangian finite element spaces defined over different triangulations of a common spatial domain, like in \autoref{fig:SketchFEUnion}. While the functions in $\Vn{1}+\dots+\Vn{N}$ are still piecewise linear, they do not always correspond to a finite element discretization on a triangulation. Still, by adding degrees of freedom one can find a triangulation and a respective linear Lagrangian finite element space $\Vplus$ containing $\Vn{1}+\dots+\Vn{N}$. 

\begin{figure}
  \begin{center}
    \includegraphics{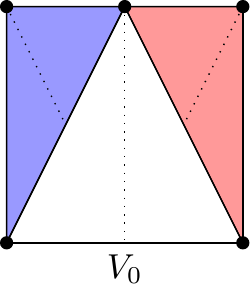}\quad
    \includegraphics{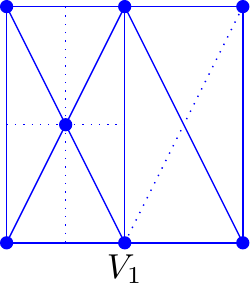}\quad
    \includegraphics{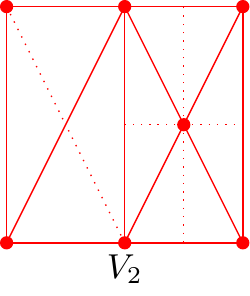}\quad
    \includegraphics{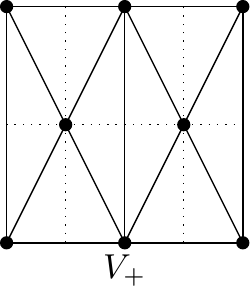}
  \end{center}
  \vspace{-5mm}
  \caption{Illustration of meshes resulting from refinement by newest vertex refinement based on a common initial triangulation corresponding to a finite element space $V_0$. Dotted lines indicate the next possible refinement. Refining the upper left triangle of the initial mesh results in $V_1$. Refining the upper right triangle of the initial mesh results in $V_2$. The common finite element space $V_+$ equals the overlay of both refined meshes, and therefore $V_+=V_1+V_2$.}
  \label{fig:SketchFEUnionNVB}
\end{figure}

A more convenient situation is encountered if the snapshots are adapted with the newest vertex bisection algorithm starting from a common initial triangulation. It is known that the smallest common refinement of two such meshes is their overlay \cite{CasconEA2008,Stevenson2007}, which implies $\Vn{1}+\dots+\Vn{N}=\Vplus$. A sketch is given in \autoref{fig:SketchFEUnionNVB}. Moreover, the mesh of $\Vplus$ can be found by repeated local refinements of any snapshot mesh. Consequently, to interpolate a function from $V_n$ to $\Vplus$ for any $n=1,\dots,N$, one only needs to interpolate the function between successive refinement steps. Because of this favorable property, we focus on refinement by newest vertex bisection in our numerical examples. However, the theory does not depend on this decision.

Besides a common finite element space of all snapshots, it will be useful to have common finite element spaces of subsets of snapshots available. To this end, we extend our notation in the following way: Let $V_{n_1},\dots,V_{n_K}$ for $K>1$ define a $K$-tuple of snapshot finite element spaces with $1\leq n_k\leq N$ for all $k=1,\dots,K$. We define $V_{n_1\dots n_K}$ as a finite element space of the same type as $V_{n_1},\dots,V_{n_K}$, with $V_{n_1}+\dots+V_{n_K}\subset V_{n_1\dots n_K}$. As a special case we denote $\Vplus = V_{1\dots N}$.

\subsection{Gramian of adapted snapshots}\label{sec:Gramian}

The first step in the computation of a POD with the method of snapshots is the creation of the snapshot Gramian, see \autoref{sec:MethodOfSnapshots}. For the case of space adapted snapshots, we consider two options: The first option is to represent all snapshots as members of a common finite element space of all snapshots. The second option is to represent pairs of snapshots as members of common finite element spaces of these pairs.

At first we provide an implementation for adaptive finite element snapshots in terms of a common finite element space of \emph{all} snapshots, where we choose $\un{1},\dots,\un{N}\in\Vplus$. We collect the finite element coefficients of the snapshots with respect to a basis of $\Vplus$ in a set of snapshot coefficient vectors $\Un{1},\dots,\Un{N}\in\mathbb R^{M_+}$ and define a snapshot matrix $\snapshotMatrix=(\Un{1},\dots,\Un{N})\in\mathbb R^{M_+\times N}$. Let $\massMatrix_+$ be the matrix associated with the $\V$-inner product of functions in $\Vplus$, so that for $u_i,u_j\in\Vplus$ we have $(u_i,u_j)_{\V}=\Un{i}^T\massMatrix_+\Un{j}$. Then the snapshot Gramian matrix is given by $\snapshotGramian = \snapshotMatrix^T\massMatrix_+\snapshotMatrix$.

Now we reformulate the computation of the snapshot Gramian so that we only need to create common finite element spaces of \emph{pairs} of snapshots. We consider the computation of a single entry of the snapshot Gramian matrix for a pair consisting of $\un{i}\in\Vn{i}$ and $\un{j}\in\Vn{j}$. Let $\Umn{ij}{i}$ and $\Umn{ij}{j}$ be the finite element coefficients of $\un{i}$ and $\un{j}$ with respect to a basis of their common finite element space $\Vn{ij}$ and let $\massMatrix_{ij}$ be the matrix associated with the $\V$-inner product of functions in $\Vn{ij}$, so that $(\un{i},\un{j})_\V=(\Umn{ij}{i})^T\massMatrix_{ij}\Umn{ij}{j}$. This means that the Gramian matrix can be filled by creating finite element spaces $\Vn{ij}$ of all pairs of snapshots. 

We have presented two ways of creating the snapshot Gramian matrix for the eigenvalue decomposition associated with the POD method of snapshots. In any case, due to the properties of the common finite element spaces we obtain the exact Gramians. The advantage of the first method is that only a single common finite element space has to be created. A possible disadvantage is that the dimension of this space may be very high. In the second method a larger number of lower-dimensional finite element spaces must be created.

\subsection{POD basis functions and approximation}

The POD basis functions are determined as linear combinations of snapshots by \eqref{eq:PODbasis}. If the snapshots are represented as members of $\Vplus$, the POD basis functions are automatically members of $\Vplus$ and can be computed by linearly combining the snapshot finite element coefficient vectors corresponding to a basis of $\Vplus$. If the snapshots are represented as members of their original adapted finite element spaces $\Vn{1},\dots,\Vn{N}$, the POD basis functions can be implicitly defined as linear combinations of snapshots. In this way, forming a basis of $\Vplus$ can be avoided, but applying a linear operator to a single POD basis function means applying this operator to all snapshots. Following this idea, there are multiple ways to represent a POD approximation: in terms of $\Vplus$, or in terms of a linear combination of POD basis functions by \eqref{eq:PODapproximation}, or in terms of snapshots by \eqref{eq:PODapproximation2}. By expressing the POD approximation and the POD basis in terms of the snapshots, one can formulate POD Galerkin models based on adaptive snapshots without the need to create the common discretization space.

The theoretical results regarding the $V$-orthogonal POD projection in the last paragraph of \autoref{sec:MethodOfSnapshots} have been stated in the context of functions in $V$. Therefore, these results do not depend on whether the snapshots have been computed with a static or an adaptive discretization. Still, the $V$-orthogonal projection requires knowledge of the function to be projected, and is therefore only valuable as a reference.

A different scenario is the computation of POD coefficients by a reduced-order model obtained via Galerkin projection. Here, knowledge of the solution is not necessary to obtain POD coefficients. As we will see in the following section, however, the snapshot discretization influences the accuracy of the POD approximation.

\section{POD Galerkin reduced-order modeling for an elliptic PDE}
\label{sec:elliptic}

We highlight the principal differences between POD Galerkin reduced-order modeling for static and adapted snapshots with an example of a parametrized elliptic boundary value problem. For the case $R=D$ we can use results from the greedy reduced basis theory \cite{AliUrbanSteih20XX}, because according to \eqref{eq:equalSpans} the POD space equals the span of the snapshots. For $R<D$ we need to take the additional POD truncation error into account. The goal of the error assessment is to understand in which way the POD truncation and the mismatch between the snapshot finite element spaces contribute to the error of the reduced solution.

\subsection{Weak formulation}

Let $\parameter\in \parameterDomain$ be a parameter vector in a domain $\parameterDomain\subset\mathbb R^K$. We define a parametrized bilinear form $a(\cdot,\cdot;\parameter)\colon\V\times\V\rightarrow\mathbb R$ which is uniformly coercive with coercivity constant 
\[
  \alpha(\parameter) = \inf_{v\in\V\setminus\{0\}}\frac{a(v,v;\parameter)}{\|v\|_V^2}\geq\alpha>0
\]
and uniformly continuous with continuity constant
\[
  \gamma(\parameter) = \sup_{v,w\in\V\setminus\{0\}}\frac{a(v,w;\parameter)}{\|v\|_V\|w\|_V}\leq\gamma<\infty.
\]
We also define a linear form $f(\cdot;\parameter)\colon\V\rightarrow\mathbb R$ which is uniformly continuous with continuity constant 
\[
  \delta(\parameter) = \sup_{v\in\V\setminus\{0\}}\frac{f(v;\parameter)}{\|v\|_V}\leq\delta<\infty.
\]

The parametrized elliptic PDE problem is now formulated as follows: For $\parameter\in\parameterDomain$, find $u(\parameter)\in\V$ such that
\begin{alignat}{2}\label{eq:ellipticPDE}
  a(u(\parameter),v;\parameter) &= f(v;\parameter)&\quad&\forall v\in\V.
\end{alignat}
The Lax-Milgram theorem guarantees well-posedness of this problem under the given continuity and coercivity assumptions. Its solution is called \emph{true solution} in the following.

\subsection{Snapshot computation}\label{sec:snapshots}

To provide snapshots for the POD computation, we introduce a collection of discretized PDE problems associated with a given discrete training set $\parameterDomain_N=\{\parameter_1,\dots,\parameter_N\}$ with $\parameter_1,\dots,\parameter_N\in\parameterDomain$. For each $\mu\in\parameterDomain_N$ we solve the respective PDE problem with an adaptive discretization scheme, which leads to the snapshot spaces $\Vn{1},\dots,\Vn{N}\subset\V$. The snapshots are solutions of the following discretized version of \eqref{eq:ellipticPDE}: For each $n=1,\dots,N$, find $\un{n}\in\Vn{n}$ such that 
\begin{alignat*}{2}
  a(\un{n},v;\parameter_n) &= f(v;\parameter_n)&\quad&\forall v\in\Vn{n}.
\end{alignat*}

For any $n=1,\dots,N$, the subspace property $\Vn{n}\subset\V$ leads to Galerkin orthogonality between the error $u(\parameter_n)-\un{n}$ and the discrete space $\Vn{n}$, which means
\begin{alignat*}{2}
  a(u(\parameter_n)-\un{n},v;\parameter_n)=0&\quad&\forall v\in\Vn{n}.
\end{alignat*}
We imply a C\'ea lemma stating
\begin{alignat}{2}\label{eq:ceaAdaptive}
  \|u(\parameter_n)-\un{n}\|_V\leq \frac{\gamma(\parameter_n)}{\alpha(\parameter_n)}\|u(\parameter_n)-P_nu(\parameter_n)\|_V,&\quad& n=1,\dots,N,
\end{alignat}
where $P_n$ denotes the $V$-orthogonal projection onto the snapshot discretization space $\Vn{n}$, so that
\[
  \|u(\parameter_n)-P_nu(\parameter_n)\|_V = \inf_{v\in\Vn{n}}\|u(\parameter_n)-v\|_V.
\]

\subsection{Reduced-order model}
\label{sec:rom}

Using the methods from \autoref{sec:Pod}, we create a POD space $\VR{R}\subset\V$ from the snapshots. The respective POD-Galerkin reduced-order model of \eqref{eq:ellipticPDE} is formulated as follows: For $\parameter\in\parameterDomain$, find $\uR{R}(\parameter)\in\VR{R}$ such that 
\begin{alignat}{2}\label{eq:rom}
  a(\uR{R}(\parameter),v;\parameter) &= f(v;\parameter)&\quad&\forall v\in\VR{R}.
\end{alignat}

The considerations regarding the computation of the snapshot Gramian in \autoref{sec:Gramian} give rise to two ways to implementing a reduced-order representation of \eqref{eq:rom}. In the first approach, the snapshots and reduced basis functions are interpreted as elements of $\Vplus$. Substituting \eqref{eq:PODapproximation} in \eqref{eq:rom} and testing against the POD basis functions leads to the following implementation: For $\parameter\in\parameterDomain$, find $\vec b(\parameter):\parameterDomain\rightarrow \mathbb R^\nPod$ such that
\begin{alignat*}{2}
  \sum_{i=1}^{R}a(\phi_i,\phi_r;\parameter)b^i(\parameter) &= f(\phi_r;\parameter)&\quad&r=1,\dots,R.
\end{alignat*}
This requires building $\Vplus$ and the respective finite element operators, which may be expensive in some cases. 

As an alternative, one can substitute \eqref{eq:PODapproximation2} in \eqref{eq:rom} and test against the POD basis functions represented in terms of the snapshots via \eqref{eq:PODbasis}, which leads to the following implementation: For $\parameter\in\parameterDomain$, find $\vec b(\parameter):\parameterDomain\rightarrow \mathbb R^\nPod$ such that
\begin{alignat*}{2}
  \sum_{i=1}^{R}\sum_{m,n=1}^N\frac{\arn{i}{m}}{\sqrt{\lambda_i}}a( \un{m}, \un{n};\parameter)\frac{\arn{r}{n}}{\sqrt{\lambda_r}}b^i(\parameter) &= \sum_{n=1}^Nf( \un{n};\parameter)\frac{\arn{r}{n}}{\sqrt{\lambda_r}}&\quad&r=1,\dots,R.
\end{alignat*}
This requires evaluating the bilinear form for all pairs of snapshots. The cost is similar to creating the snapshot Gramian needed for the POD computation.

In any case, further information about the dependence of the linear and bilinear forms on the parameter $\parameter$ is needed to obtain a reduced-order model which can be evaluated for any $\parameter\in\parameterDomain$ at a cost which does not depend on the number of spatial degrees of freedom. An example is given in \autoref{sec:ConvectionDiffusion}.

\subsection{Error assessment}

For the error assessment of $\uR{R}$ from \eqref{eq:rom}, we first recall the main results for static discretizations. Then we study the adaptive case and point out major differences compared to the static case. We restrict our attention to $\parameter\in\parameterDomain_N$. For practical applications this means that a sufficiently rich snapshot set is assumed, so that the error from discretizing $\parameterDomain$ is negligible.

\subsubsection{Static discretization}
\label{sec:errorStatic}

Assume $\VR{R}\subset\Vn{n}$ for $n=1,\dots,N$. This assumption holds if the snapshots have been computed with static finite elements. 

At first we study the error between the snapshots and the reduced-order solution evaluated at the corresponding training parameter values. From $\VR{R}\subset\Vn{n}$, we can derive a Galerkin orthogonality between this error and the POD space,
\begin{alignat*}{2}
  a(\un{n}-\uR{R}(\parameter_n),v;\parameter_n)=0&\quad&\forall v\in\VR{R},\quad n=1,\dots,N.
\end{alignat*}
The C\'ea lemma, following from coercivity and continuity, states the relation to the POD approximation error,
\begin{alignat*}{2}
  \|\un{n}-\uR{R}(\parameter_n)\|_V\leq \frac{\gamma(\parameter_n)}{\alpha(\parameter_n)}\|\un{n}-P^R\un{n}\|_V,&\quad&n=1,\dots,N.
\end{alignat*}
Moreover, \eqref{eq:sumEig} implies $\uR{D}(\parameter_n) = \un{n}$ for $n=1,\dots,N$. The fact that the snapshots are recovered for large enough $R$ is called asymptotic snapshot reproducibility.

For the error with respect to the true solution we use a triangle inequality and respective C\'ea lemmas to obtain
\begin{alignat*}{2}
  \|u(\parameter_n)-\uR{R}(\parameter_n)\|_V&\leq \|u(\parameter_n)-\un{n}\|_\V+\|\un{n}-\uR{R}(\parameter_n)\|_V\\
                              &\leq \frac{\gamma(\parameter_n)}{\alpha(\parameter_n)}\|u(\parameter_n)-P_nu(\parameter_n)\|_V+\frac{\gamma(\parameter_n)}{\alpha(\parameter_n)}\|\un{n}-P^R\un{n}\|_V
\end{alignat*}
for $n=1,\dots,N$. For maximum efficiency, the errors of the finite element discretization and the POD truncation should be balanced. It should be noticed, however, that the first term stems from the off-line discretization, which is only relevant for the setup of the reduced-order model, while the second term stems from the on-line discretization, which is also relevant for the evaluation time of the reduced-order model.

\subsubsection{Adaptive discretization}
\label{sec:errorAdaptive}

We derive error inequalities similar to the ones in \autoref{sec:errorStatic}, but for the more general case, where the snapshots are members of different finite element spaces. The assumption $\VR{R}\subset\Vn{n}$ for $n=1,\dots,N$ is usually not satisfied in the adaptive case, we only have $\VR{R}\subset\Vplus$ and $\Vn{n}\subset\Vplus$. As a consequence, we are not able to use a Galerkin orthogonality between the reduced-order error $\un{n}-\uR{R}(\parameter_n)$ and the reduced space $\VR{R}$ for $n=1,\dots,N$. 

We start with the error between the solution of the reduced-order model and the true solution. Due to $\VR{R}\subset\V$, for any $\parameter\in\parameterDomain$ one can derive a Galerkin orthogonality
\begin{alignat*}{2}
  a(u(\parameter)-\uR{R}(\parameter),v;\parameter)=0&\quad&\forall v\in\VR{R}
\end{alignat*}
and a corresponding C\'ea lemma
\begin{alignat*}{2}
  \|u(\parameter)-\uR{R}(\parameter)\|_V\leq \frac{\gamma(\parameter)}{\alpha(\parameter)}\|u(\parameter)-P^Ru(\parameter)\|_V.
\end{alignat*}
We split the right-hand side of the C\'ea lemma into contributions from the snapshot computation and from the POD truncation. To exclude the error associated with the discretization of the parameter domain, we consider only $\parameter\in\parameterDomain_N$. The derivation starts with adding a zero to the right-hand side of the C\'ea lemma for the reduced-order model and subsequently uses triangle inequalities and the properties of orthogonal projections,
\begin{alignat*}{1}
  \|u(\parameter_n)-\uR{R}(\parameter_n)\|_V&\leq \frac{\gamma(\parameter_n)}{\alpha(\parameter_n)}\|u(\parameter_n)-u_n+u_n-P^Ru_n+P^Ru_n-P^Ru(\parameter_n)\|_V\\
  &\leq \frac{\gamma(\parameter_n)}{\alpha(\parameter_n)}\|u_n-P^Ru_n\|_V
  +    \frac{\gamma(\parameter_n)}{\alpha(\parameter_n)}\|(I-P^R)(u(\parameter_n)-u_n)\|_V\\
  &\leq \frac{\gamma(\parameter_n)}{\alpha(\parameter_n)}\|u_n-P^Ru_n\|_V
  +    \frac{\gamma(\parameter_n)}{\alpha(\parameter_n)}\|u(\parameter_n)-u_n\|_V
\end{alignat*}
for $n=1,\dots,N$. Using \eqref{eq:ceaAdaptive} we obtain
\begin{alignat*}{1}
  \|u(\parameter_n)-\uR{R}(\parameter_n)\|_V&\leq \frac{\gamma(\parameter_n)}{\alpha(\parameter_n)}\|u_n-P^Ru_n\|_V + \frac{\gamma(\parameter_n)^2}{\alpha(\parameter_n)^2}\|u(\parameter_n)-P_nu(\parameter_n)\|_V
\end{alignat*}
for $n=1,\dots,N$. This means that for parameter values in $\parameterDomain_N$, the error between the true and the reduced-order solution can be split into contributions from the projection of the respective snapshot onto the POD space and from the projection of the true solution onto the respective snapshot finite element space. In absence of POD truncation, i.e. for $R=D$, the POD projection error vanishes and we obtain a variant of known results from greedy reduced basis theory \cite{AliUrbanSteih20XX}.

Because in general $\VR{R}\not\subset\Vn{n}$ for $n=1,\dots,N$, we are not able to derive a C\'ea lemma for $\un{n}-\uR{R}(\parameter_n)$. A straight-forward approach is using the results from above to obtain
\begin{align*}
  &\|\un{n}-\uR{R}(\parameter_n)\|_V\\
  &\qquad\leq\|\un{n}-u(\parameter_n)\|_V+\|u(\parameter_n)-\uR{R}(\parameter_n)\|_V\\  
  &\qquad\leq\frac{\gamma(\parameter_n)}{\alpha(\parameter_n)}\left(1+\frac{\gamma(\parameter_n)}{\alpha(\parameter_n)}\right)\|u(\parameter_n)-P_nu(\parameter_n)\|_V+\frac{\gamma(\parameter_n)}{\alpha(\parameter_n)}\|u_n-P^Ru_n\|_V.
\end{align*}
Alternatively, we can use another result adapted from the literature \cite{AliUrbanSteih20XX}: Due to coercivity, continuity, $u_n\in\VR{D}$ and Galerkin orthogonality between $\VR{D}$ and the error between the true solution and the solution of the reduced-order model we have
\begin{alignat*}{1}
  \alpha(\parameter_n)\|\un{n}-\uR{D}(\parameter_n)\|_V^2&\leq a(\un{n}-\uR{D}(\parameter_n),\un{n}-\uR{D}(\parameter_n);\parameter_n)\\
  &= a(\un{n}-u(\parameter_n),\un{n}-\uR{D}(\parameter_n);\parameter)\\
  &\leq \gamma(\parameter_n)\|\un{n}-u(\parameter_n)\|_V\|\un{n}-\uR{D}(\parameter_n)\|_V,\quad n=1,\dots,N,
\end{alignat*}
so that
\begin{alignat*}{1}
  \|\un{n}-\uR{D}(\parameter_n)\|_V&\leq \frac{\gamma(\parameter_n)}{\alpha(\parameter_n)}\|\un{n}-u(\parameter_n)\|_V,\quad n=1,\dots,N.
\end{alignat*}
Because $\VR{R}\subset\VR{D}$ for $R\leq D$ we obtain
\begin{align*}
  &\|\un{n}-\uR{R}(\parameter_n)\|_V\\
  &\qquad\leq\|\un{n}-\uR{D}(\parameter_n)\|_V+\|\uR{D}(\parameter_n)-\uR{R}(\parameter_n)\|_V\\  
  &\qquad\leq \frac{\gamma(\parameter_n)}{\alpha(\parameter_n)}\|\un{n}-u(\parameter_n)\|_V +\frac{\gamma(\parameter_n)}{\alpha(\parameter_n)}\|\uR{D}(\parameter_n)-P^R\uR{D}(\parameter_n)\|_V\\
  &\qquad\leq\frac{\gamma(\parameter_n)^2}{\alpha(\parameter_n)^2}\|u(\parameter_n)-P_nu(\parameter_n)\|_V+\frac{\gamma(\parameter_n)}{\alpha(\parameter_n)}\|\uR{D}(\parameter_n)-P^R\uR{D}(\parameter_n)\|_V,
\end{align*}
with the help of the C\'ea lemma of the snapshot computation. In any case, the error between a snapshot and the solution of the reduced-order model at the corresponding parameter value contains a component from the finite element computation, which was not present for static snapshots.

\section{Numerical example of a convection-diffusion equation}
\label{sec:ConvectionDiffusion}

We apply POD model order reduction to a two-dimensional convection-diffusion problem, where the transport direction serves as a parameter and adaptive finite element snapshots are taken over the parameter interval. The test case illustrates the computational methods introduced in \autoref{sec:Pod} and the theoretical results derived in \autoref{sec:elliptic}.

\subsection{Problem setting}\label{sec:ConvectionDiffusionProblem}

We consider a pa\-ra\-me\-tri\-zed boundary value problem based on a convection-diffusion equation in two dimensions,
\begin{alignat*}{2}
  v_x(\mu)\partial_x u + v_y(\mu)\partial_y u - \nu\partial_{xx} u - \nu\partial_{yy} u &= 1 &\quad&\forall (\vec x,\parameter)\in \Omega\times\parameterDomain,
\end{alignat*}
with solution $u(\vec x,\parameter):\Omega\times\parameterDomain\rightarrow \mathbb R$ for a spatial domain $\Omega=(0,1)\times(0,1)$, a parameter interval $\parameterDomain=[0,1]$ and a diffusivity of $\nu=0.01$. The dependence of the convective velocity components $v_x(\mu)$ and $v_y(\mu)$ on the parameter $\mu$ is given by $v_x = \cos(0.25\pi\mu)$ and $v_y = \sin(0.25\pi\mu)$, which means that only the direction of the velocity vector is varied. At the boundary $\partial\Omega$ of the spatial domain $\Omega$, we specify homogeneous Dirichlet conditions, so that
\begin{alignat*}{2}
  u(\vec x,\parameter) &= 0 &\quad&\forall (\vec x,\parameter)\in \partial\Omega\times\parameterDomain.
\end{alignat*}
A spatial weak form of the parametrized convection-diffusion problem is given by \eqref{eq:ellipticPDE}, where  $V:=H^1_0(\Omega)$ and the bilinear and linear forms are defined by
\begin{align*}
  a(w,\psi;\mu)   &:= \int_\Omega v_x(\mu)\partial_x w\,\psi + v_y(\mu)\partial_y w\,\psi + \nu\partial_x w\,\partial_x \psi + \nu\partial_y w\,\partial_y \psi\; \mathrm d\vec x,\\
  f(\psi) &:= \int_\Omega \psi\; \mathrm d\vec x.
\end{align*}

\subsection{Discretization}\label{sec:convectionDiffusionDiscretization}

The general procedure for the snapshot computation and the discretization of the parameter domain has been described in \autoref{sec:snapshots}. Here, we focus on the weak form of the convection-diffusion equation. To provide a complete description of the snapshot computation, the choice of $\Vn{1},\dots,\Vn{\nSnap}$ is characterized in the following.

For the numerical discretization, we use a custom implementation \cite{Ullmann2016} of piecewise linear Lagrangian finite elements on triangular meshes. As an error indicator, we employ the Matlab function \texttt{pdejmps} \cite{pdetool}, defining an error measure for a triangle $K$ by
\[
  E(K) = \left(\frac12\sum_{\tau\in\partial K}h_\tau^2[\vec n_\tau\cdot(\nu\nabla u_h)]^2\right)^\frac12,
\]
where $\tau$ denotes an edge along the triangle boundary $\partial K$, $h_\tau$ is the length of the edge, $\vec n_\tau$ is the unit normal, and $[\cdot]$ denotes a jump across the edge. The refinement loop is started from a coarse initial mesh which is identical for all snapshot simulations. In each step, the triangles with the largest error contributions are refined by newest vertex bisection. The marking of triangles and the termination of the refinement loop is done with Matlab \texttt{pdeadgsc} \cite{pdetool}.

\begin{figure}
    \includegraphics[scale = 0.3,trim=80 20 120 0,clip = true]{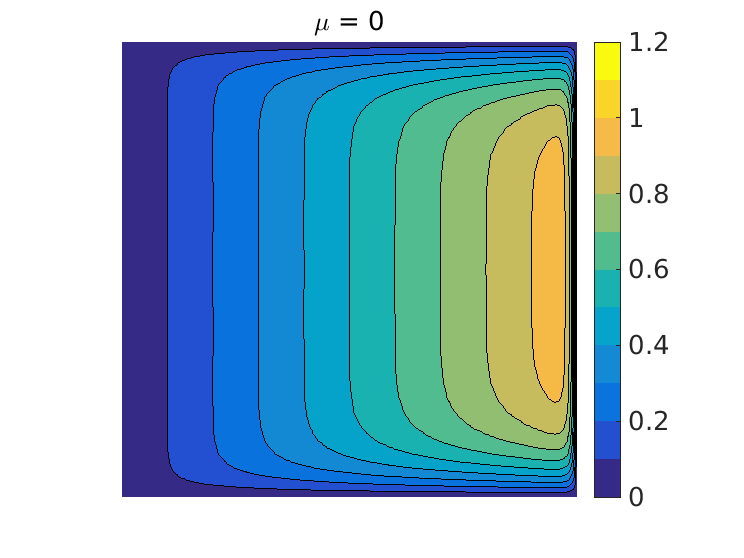}
    \includegraphics[scale = 0.3,trim=80 20 120 0,clip = true]{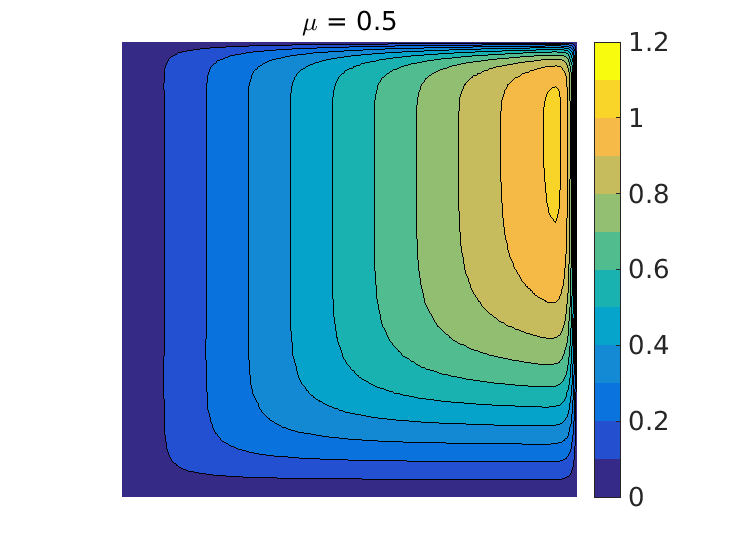}
    \includegraphics[scale = 0.3,trim=80 20  50 0,clip = true]{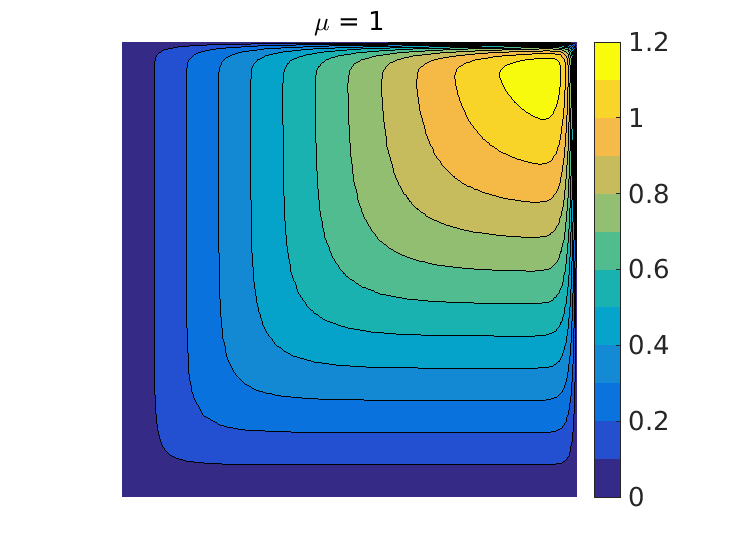}
  \caption{Solution of the convection-diffusion problem for different parameter values.}
  \label{fig:solutionConvectionDiffusion}
\end{figure}

\begin{figure}
    \includegraphics[scale = 0.3,trim=80 20 120 0,clip = true]{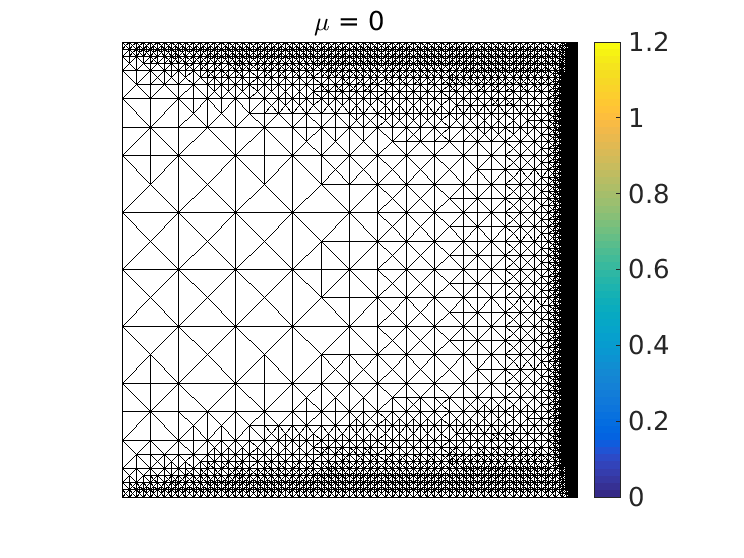}
    \includegraphics[scale = 0.3,trim=80 20 120 0,clip = true]{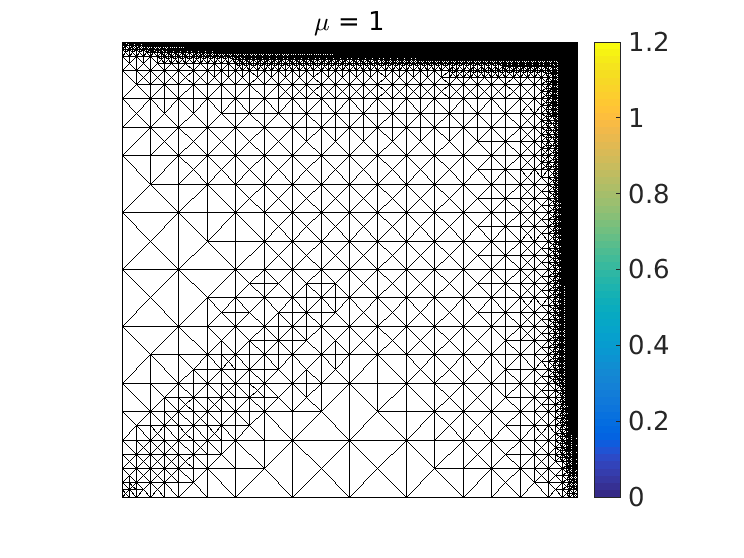}
    \includegraphics[scale = 0.3,trim=80 20 120 0,clip = true]{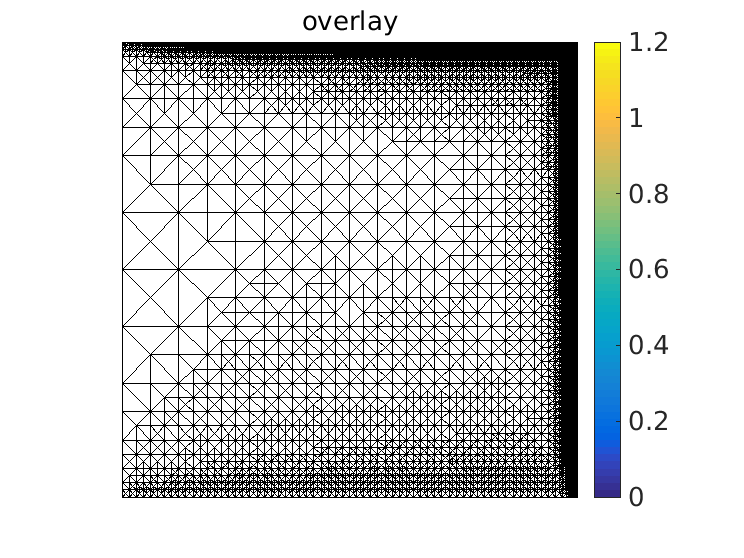}\\
  \caption{Adapted triangulations corresponding to \autoref{fig:solutionConvectionDiffusion}, and the overlay of grids of all snapshots, consituting $\Vplus$.}
  \label{fig:meshesConvectionDiffusion}
\end{figure}

Numerical solutions for varying parameter values are presented in \autoref{fig:solutionConvectionDiffusion}. Adapted spatial meshes are plotted in \autoref{fig:meshesConvectionDiffusion} together with an overlay of all snapshot meshes constituting $\Vplus$ for further reduced-order modeling.

\subsection{Reduced-order modeling}

In order to generate snapshots for subsequent reduced-order modeling, we solve the problem for 33 parameter values distributed equidistantly over the parameter interval and store the respective solutions. We compute a POD of these snapshots and create Galerkin reduced-order models of varying dimension, using the same techniques as described in \autoref{sec:rom}.
We can rewrite the model as an equation for the POD coefficient vector: Find $\vec b(\parameter):\parameterDomain\rightarrow \mathbb R^\nPod$ such that
\begin{alignat*}{2}
  (A_x^Rv_x(\mu) + A_y^Rv_y(\mu) + A^R_\nu) \vec b &= \vec F^R.
\end{alignat*}
Expressions for the constant model coefficient matrices $A_x^R$, $A_y^R$, $A_\nu^R$ and the right-hand side vector $\vec F^R$ follow from substitution of the POD expansion \eqref{eq:PODapproximation} into \eqref{eq:rom} and testing against the POD basis functions. In this example, refinement was mostly necessary near the boundary. Therefore, forming the common finite element space of all snapshots was not very costly. Consequently, the snapshot Gramian and the coefficients of the reduced-order model were created using $\Vplus$.

\subsection{Results}\label{sec:convectionDiffusionResults}

In order to test how well the reduced-order models reproduce the underlying snapshots, we solve the models for the parameter values in the training set. The error between the snapshots and the solution of the POD-Galerkin model is named $\epsilon_\text{ROM}$. It can be viewed as a benchmark for the Galerkin reduced-order model under the assumption that sufficient snapshot data has been used to create the POD. It is measured in the relative norm induced by the POD, so that 
\begin{equation}\label{eq:e_ROM}
  \epsilon_\text{ROM} = \sqrt{\sum_{\iSnap=1}^{\nSnap}\|u_{\iSnap} - \uR{\nPod}(\parameter_\iSnap)\|_V^2}\bigg/\sqrt{\sum_{\iSnap=1}^{\nSnap}\|u_{\iSnap}\|_V^2}.
\end{equation}

The error between the snapshots and their orthogonal projection on the POD is named $\epsilon_\text{POD}$. It can be viewed as a benchmark for the POD in the sense that it characterizes the ability of the snapshot data to be represented in a low-dimensional space. It is measured in the relative norm induced by the POD, which enables an alternative expression in terms of the POD eigenvalues via \eqref{eq:sumEig}:
\begin{equation}\label{eq:e_POD}
  \epsilon_\text{POD} = \sqrt{\sum_{\iSnap=1}^{\nSnap}\|u_{\iSnap} - P^Ru_{\iSnap}\|_V^2}\bigg/\sqrt{\sum_{\iSnap=1}^{\nSnap}\|u_{\iSnap}\|_V^2} = \sqrt{\sum_{n=R+1}^D\lambda_n}\bigg/\sqrt{\sum_{n=1}^D\lambda_n}.
\end{equation}

To compare the approximation errors resulting from the model order reduction with the finite element discretization error, we computed reference snapshots $u_1^\text{ref},\dots,u_N^\text{ref}$ with a stricter spatial tolerance. The resulting estimated finite element discretization error of the original snapshots is measured in the same norm as above, so that
\begin{equation}\label{eq:e_FEM}
  \epsilon_\text{FEM} = \sqrt{\sum_{\iSnap=1}^{\nSnap}\|u_{\iSnap}^\text{ref} - u_{\iSnap}\|_V^2}\bigg/\sqrt{\sum_{\iSnap=1}^{\nSnap}\|u_{\iSnap}^\text{ref}\|_V^2}.
\end{equation}

\begin{figure}
  \begin{center}
    \includegraphics[scale = 0.3,trim=10 0 10 0,clip = false]{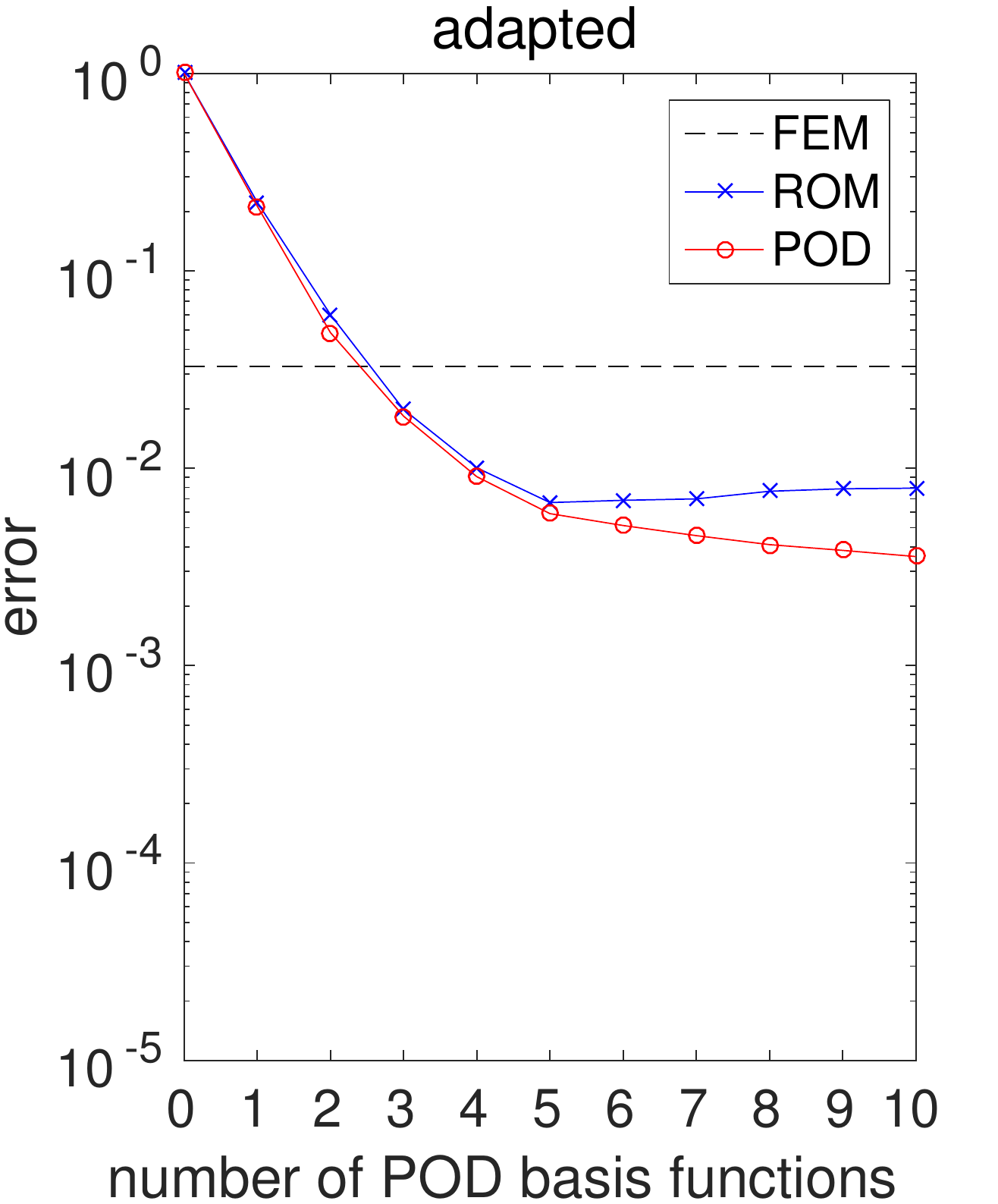}
    \quad
    \includegraphics[scale = 0.3,trim=20 0 10 0,clip = true]{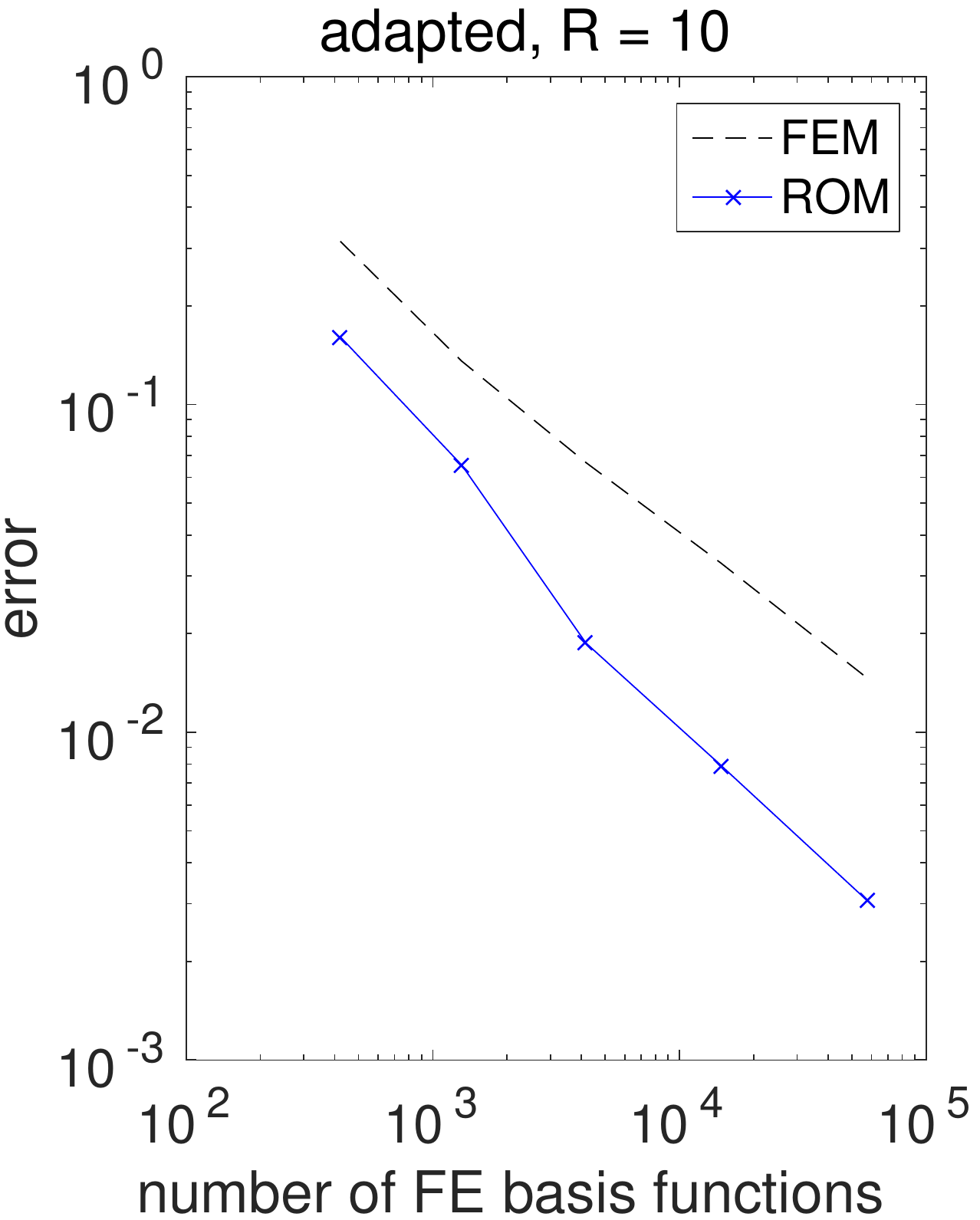}
    \quad
    \includegraphics[scale = 0.3,trim=20 0 10 0,clip = true]{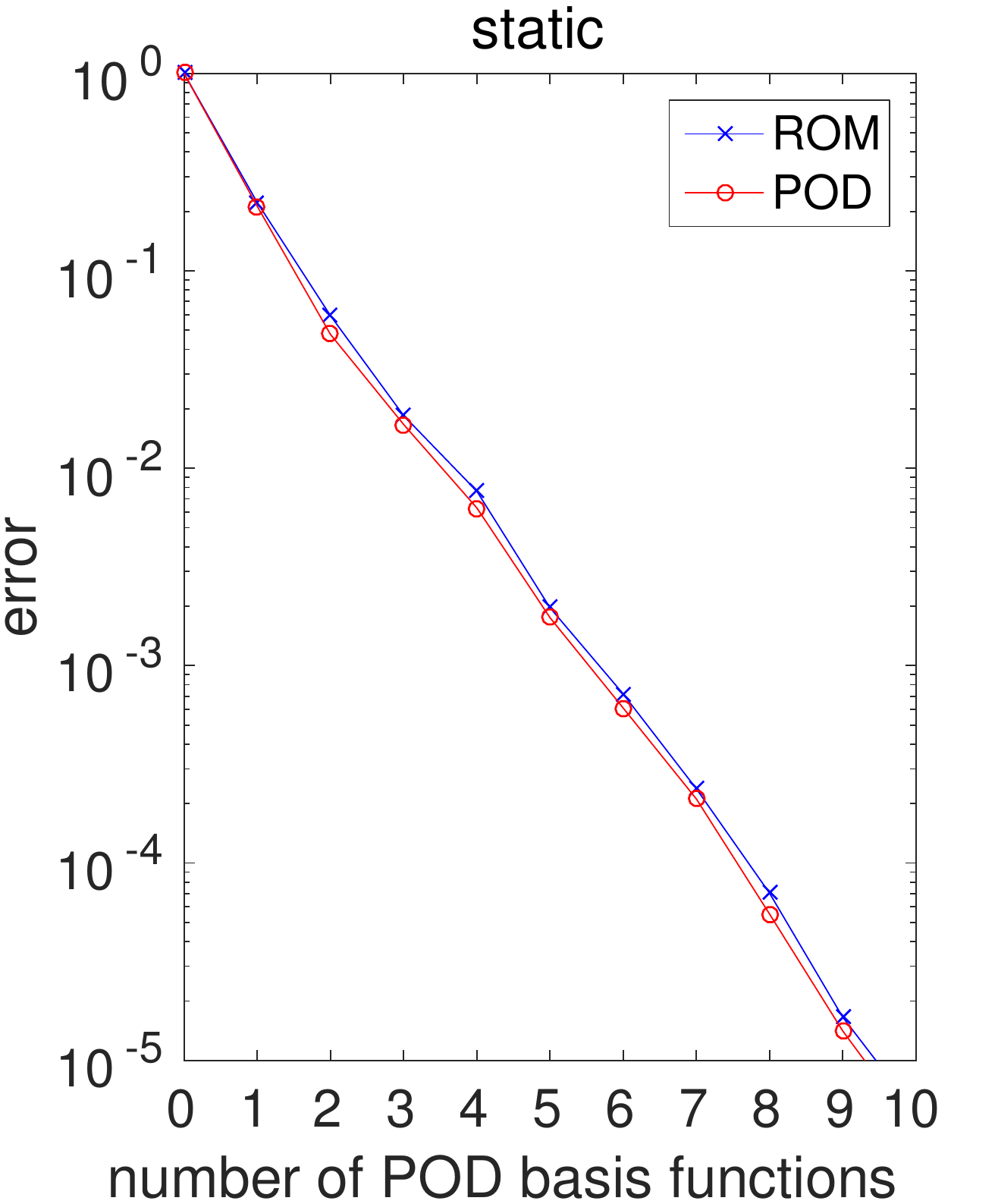}
  \end{center}
  \caption{Error $\epsilon_\text{POD}$ of the POD projection, error $\epsilon_\text{ROM}$ of the solution of the reduced-order model, and finite element discretization error $\epsilon_\text{FEM}$, corresponding to \eqref{eq:e_ROM}--\eqref{eq:e_FEM}. Left: adapted snapshots, middle: dependence on the spatial refinement for 10 POD basis functions, right: static snapshots computed on the common grid of the space-adapted snapshots.}
  \label{fig:convergenceConvectionDiffusion}
\end{figure}

The errors defined in \eqref{eq:e_ROM}--\eqref{eq:e_FEM} are plotted in \autoref{fig:convergenceConvectionDiffusion} (left). We observe that the POD projection error $\epsilon_\text{POD}$ decreases monotonically, which is implied by its representation in terms of the POD eigenvalues. The convergence of the error of the reduced-order model $\epsilon_\text{ROM}$, however, stagnates at some level.

\autoref{fig:convergenceConvectionDiffusion} (middle) shows that the stagnation level is related to the finite element error of the snapshots, as suggested by the results of section \ref{sec:errorAdaptive}. In the plot, the number of POD basis functions is kept fixed, but the number of finite element functions is varied by changing the tolerance of the refinement algorithm. One can observe that the error of the reduced-order model and the POD projection error roughly follow the finite element discretization error up to some constant.

The error depending on to the number of basis functions is shown in \autoref{fig:convergenceConvectionDiffusion} (right) for simulations on a fixed grid given by the overlay of all adapted snapshots. Here both the error of the POD projection and the solution of the reduced-order model converge. This can be explained by the Galerkin orthogonality property detailed in section \ref{sec:errorStatic}. The convergence rate is higher in the case of a static finite element discretization. Our interpretation of this behavior is that in the adaptive case the POD basis functions of higher index start approximating spatial artifacts resulting from the different discretizations. 

\begin{figure}
    \includegraphics[scale = 0.3,trim=30 10 90 5,clip = true]{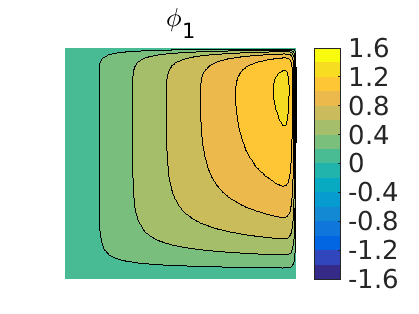}
    \includegraphics[scale = 0.3,trim=30 10 90 5,clip = true]{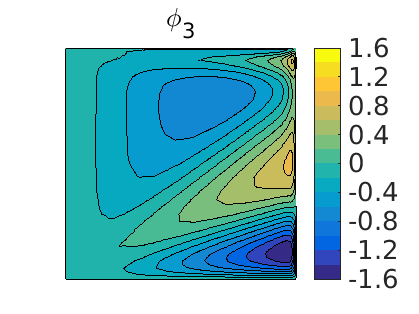}
    \includegraphics[scale = 0.3,trim=30 10 90 5,clip = true]{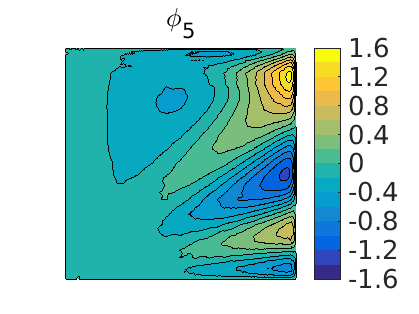}
    \includegraphics[scale = 0.3,trim=30 10 90 5,clip = true]{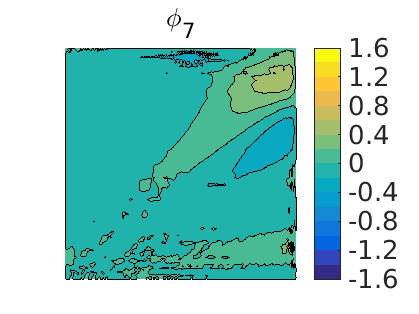}
    \includegraphics[scale = 0.3,trim=30 10  0 5,clip = true]{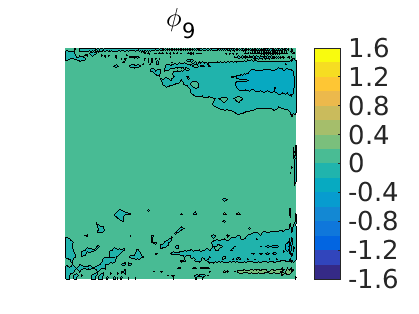}\\[3mm]
    \includegraphics[scale = 0.3,trim=30 10 90 5,clip = true]{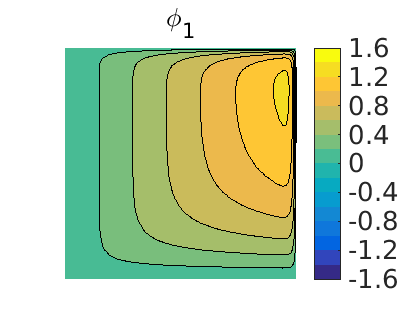}
    \includegraphics[scale = 0.3,trim=30 10 90 5,clip = true]{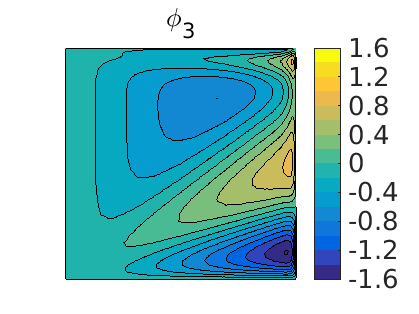}
    \includegraphics[scale = 0.3,trim=30 10 90 5,clip = true]{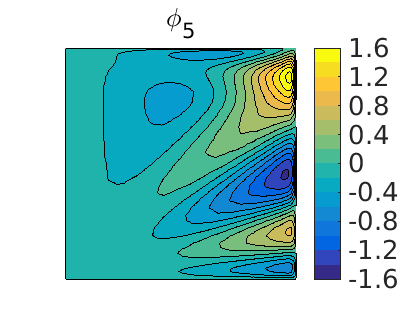}
    \includegraphics[scale = 0.3,trim=30 10 90 5,clip = true]{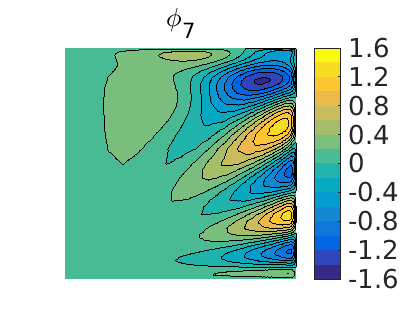}
    \includegraphics[scale = 0.3,trim=30 10  0 5,clip = true]{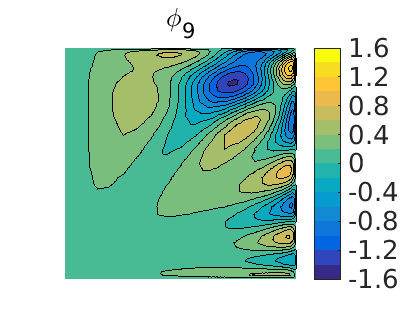}
  \caption{A selection of POD basis functions resulting from 33 snapshots. Top: adapted, bottom: static.}
  \label{fig:basisConvectionDiffusion}
\end{figure}

To underline our interpretation of the differences in the POD approximation error $\epsilon_\text{POD}$ between adaptive and static snapshots simulations, we plot a selection of POD basis functions resulting from the adapted snapshots in \autoref{fig:basisConvectionDiffusion} (top). We observe that the POD basis functions corresponding to the adaptive simulation start exhibiting local variations in the size of typical mesh cells when the index is increased. For comparison, we present a selection of POD basis functions resulting from the static snapshots in \autoref{fig:basisConvectionDiffusion} (bottom). It is reasonable that the oscillations visible in these plots are necessary to approximate the parameter-dependent physical structures of the solution with increasing accuracy, see \autoref{fig:solutionConvectionDiffusion}. 

Despite this qualitative differences in the appearance of the POD basis functions of higher index, we stress that in this example including more than 4 basis functions does not decrease the total error of the reduced-order solution. This is because the POD errors are dominated by the finite element approximation error for $R\geq 4$ in both the adaptive and the static case, see \autoref{fig:convergenceConvectionDiffusion} (left).

\section{Numerical example of a Burgers equation}
\label{sec:Burgers}

To illustrate the potential of POD-Galerkin modeling with adaptive snapshots for non-linear time-dependent problems, we apply our techniques to a Burgers equation. Adaptive finite element snapshots are taken over the time and parameter domain. Note that by choosing the parameter domain equal to a single point, a time-dependent non-parametrized Burgers problem is obtained as a special case. The test case suggests that the theoretical results derived in \autoref{sec:elliptic} for a linear elliptic setting can be transferred to a parabolic setting.

\subsection{Problem setting}

We consider a pa\-ra\-me\-tri\-zed initial-boundary value problem based on a scalar viscous Burgers equation in two dimensions,
\begin{alignat*}{2}
  \partial_t u + u\partial_x u - \nu\partial_{xx} u - \nu\partial_{yy} u &= 0 &\quad&\forall (t,\vec x,\parameter)\in I\times \Omega\times\parameterDomain,
\end{alignat*}
with solution $u(t,\vec x,\parameter):I\times\Omega\times\parameterDomain\rightarrow \mathbb R$ for a time interval $I = (0,1.2]$, a spatial domain $\Omega=(0,1)\times(0,0.5)$, and a parameter interval $\parameterDomain = [0,1]$. The problem is a modified version of a numerical example in \cite{DrohmannEA2012}. A parametrized initial condition is given by
\begin{align*}
  u(0,\vec x,\parameter) - u_0(\vec x)\parameter = 0\quad\forall (\vec x,\parameter)\in\Omega\times\parameterDomain
\end{align*}
with spatial data $u_0(x):\Omega\rightarrow \mathbb R$ given by a two-dimensional sinusoidal profile
\[
  u_0 = 0.5+0.5\sin\big((x-y-0.75)\pi\big)\sin\big((x+y+0.25)\pi\big).
\]
At the boundary we choose a combination of homogeneous Neumann and periodic conditions, so that
\begin{alignat*}{2}
  u(t,(0,y)^T,\parameter) - u(t,(1,y)^T,\parameter)&= 0 &\quad&\forall (t,y,\parameter)\in I\times (0,0.5)\times\parameterDomain,\\
  \partial_y u(t,(x,0)^T,\parameter) = \partial_y u(t,(x,0.5)^T,\parameter) &= 0&&\forall(t,x,\parameter)\in I\times (0,1)\times\parameterDomain.
\end{alignat*}
We choose a viscosity constant of $\nu=0.001$. Note that for the inviscid Burgers equation with $\nu=0$ we would have $u(t,\vec x,\mu) = \mu u(t\mu,\vec x,1)$.

In the context of \autoref{sec:Burgers}, we define $V\subset H^1(\Omega)$ as the Sobolev space of $L^2(\Omega)$ functions which have weak first derivatives in $L^2(\Omega)$ and which satisfy the periodic boundary conditions in a weak sense. A spatial weak form of the parametrized Burgers problem is then given as follows: For $u_0\in L^2(\Omega)$, find $u(t,\parameter):I\times\parameterDomain\rightarrow V$, such that
\begin{alignat*}{3}
  (u(0,\parameter)-u_0\parameter,\phi)&= 0 &&\forall\phi\in V,&&\parameter\in\parameterDomain,\\
  (\partial_t u,\phi) + a(u,\phi) + b(u,u,\phi) &= 0 &\quad&\forall\phi\in V,&\quad&(t,\parameter)\in I\times\parameterDomain,
\end{alignat*}
where $(\cdot,\cdot)$ denotes the $L^2(\Omega)$-inner product and the bilinear and trilinear forms are defined by
\begin{align*}
  a(v,\psi)   &:= \nu\int_\Omega \partial_x v\,\partial_x \psi + \partial_y v\,\partial_y \psi\; \mathrm d\vec x,&
  b(v,w,\psi) &:= \int_\Omega v\partial_x w\,\psi\; \mathrm d\vec x.
\end{align*}

\subsection{Discretization}

We apply Rothe's method to the problem, i.e.\ we employ a time discretization followed by space discretizations of the sequence of resulting boundary value problems. We use a constant time step size $\tau = T/(\nTime-1)$ with according discrete time instances $t_\iTime = (\iTime-1)\tau$ and semi-discretized solutions $\tilde u_{\iTime}\approx u(t_\iTime,\cdot)$ for $\iTime=1,\dots,\nTime$. An implicit Euler time discretization is given as follows: For $u_0\in L^2(\Omega)$, find $\tilde u_{1}(\parameter),\dots,\tilde u_{\nTime}(\parameter):\parameterDomain\rightarrow V$ such that
\begin{alignat}{3}
  (\tilde u_{\iTime}-u_0\parameter,\phi)&= 0 &\quad&\forall\phi\in V,&\quad&\iTime=1,\label{eq:semiTime1}\\
  (\tilde u_{\iTime}-\tilde u_{\iTime-1},\phi) + \tau a(\tilde u_{\iTime},\phi) + \tau b(\tilde u_{\iTime},\tilde u_{\iTime},\phi) &= 0 &\quad&\forall\phi\in V,&&\iTime=2,\dots,\nTime.\label{eq:semiTime2}
\end{alignat}
This semi-discretization is the starting point for further adaptive space discretization and reduced-order modeling.

In order to obtain snapshots for a reduced-order model, we discretize the parameter domain into a training set $\parameterDomain_\nParam:=\{\parameter_1,\dots,\parameter_{\nParam}\}$ with $\parameter_1,\dots,\parameter_{\nParam}\in\parameterDomain$. We introduce a discretization in space with adaptive finite elements, so that each snapshot belongs to an individual finite element space chosen with respect to some global error tolerance. We denote the respective fully discrete solutions and finite element spaces by $u_{\iTime,\iParam}\in V_{\iTime,\iParam}$ for $\iTime=1,\dots,\nTime$ and $\iParam=1,\dots,\nParam$, so that $u_{\iTime,\iParam}\approx u(t_\iTime,\parameter_\iParam)$. A discretization of the problem with an implicit Euler scheme in time and adaptive finite elements in space for the training parameter values is given as follows: For $u_0\in L^2(\Omega)$ and $\iParam=1,\dots,\nParam$, find $u_{1,\iParam}\in V_{1,\iParam},\dots,u_{\nTime,\iParam}\in V_{\nTime,\iParam}$ such that
\begin{alignat*}{3}
  (u_{k,\iParam}-u_0\parameter_\iParam,\phi)&= 0 &\;&\forall\phi\in V_{k,\iParam},&\;&\iTime=1,\\
  (u_{\iTime,\iParam}-u_{\iTime-1,\iParam},\phi) + \tau a(u_{\iTime,\iParam},\phi) + \tau b(u_{\iTime,\iParam},u_{\iTime,\iParam},\phi) &= 0 &\quad&\forall\phi\in V_{\iTime,\iParam},&&\iTime=2,\dots,\nTime.
\end{alignat*}
The finite element discretization leads to a set of non-linear algebraic equations which can be solved in parallel for $\iParam=1,\dots,\nParam$ and in sequence for $\iTime=1,\dots,\nTime$. The non-linear equations are solved with a standard Newton method using a direct solver for the resulting sequence of linear algebraic equations in each time step.

\begin{figure}
  \includegraphics[scale = 0.3,trim=85 25 110 20,clip = true]{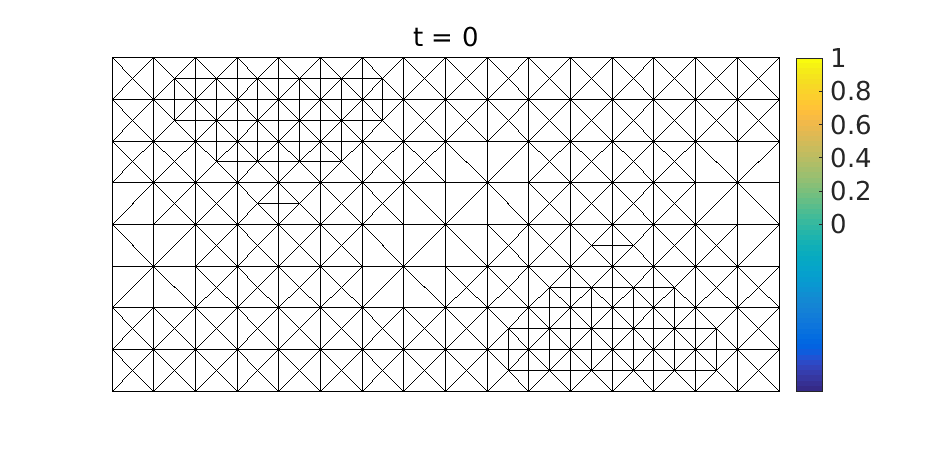}\quad
  \includegraphics[scale = 0.3,trim=85 25  45 20,clip = true]{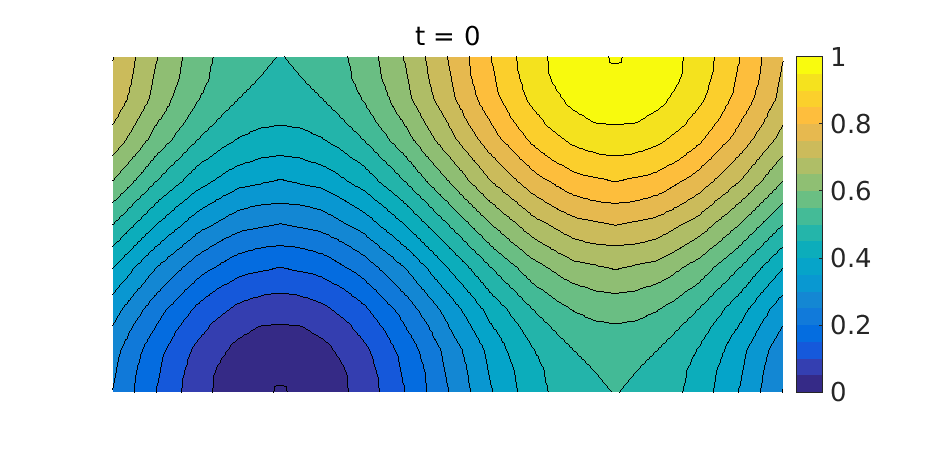}\\
  \includegraphics[scale = 0.3,trim=85 25 110 20,clip = true]{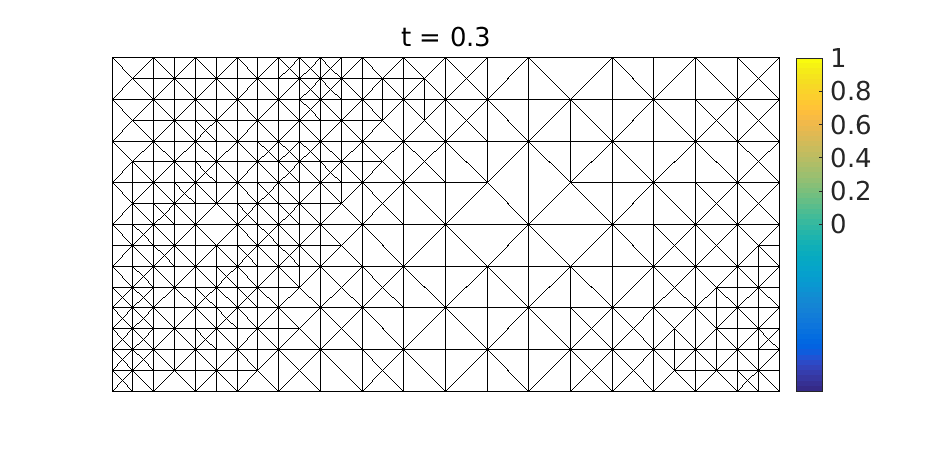}\quad
  \includegraphics[scale = 0.3,trim=85 25  45 20,clip = true]{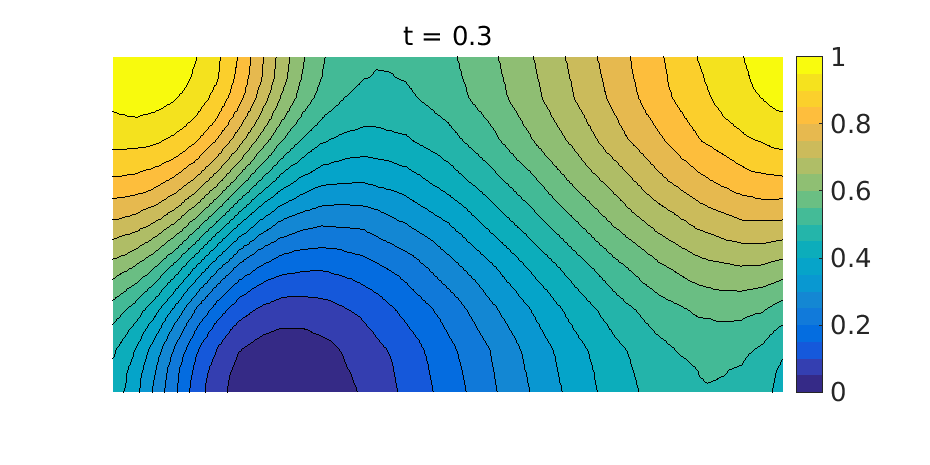}\\  \includegraphics[scale = 0.3,trim=85 25 110 20,clip = true]{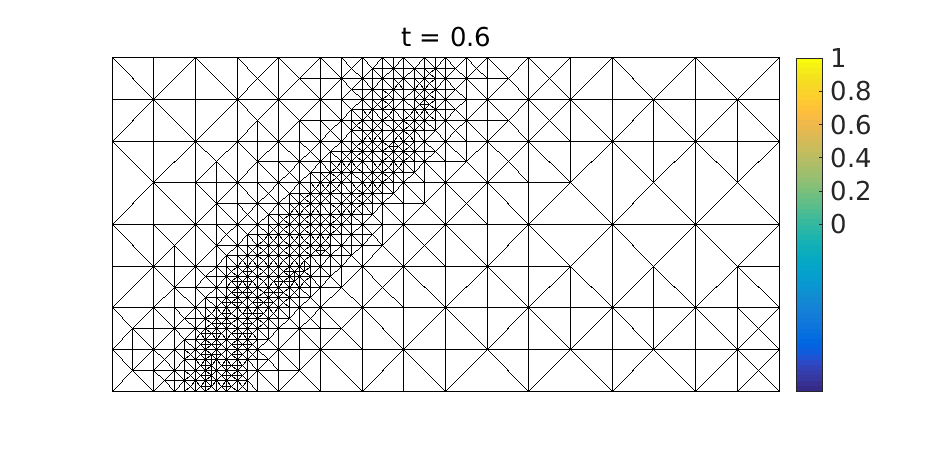}\quad
  \includegraphics[scale = 0.3,trim=85 25  45 20,clip = true]{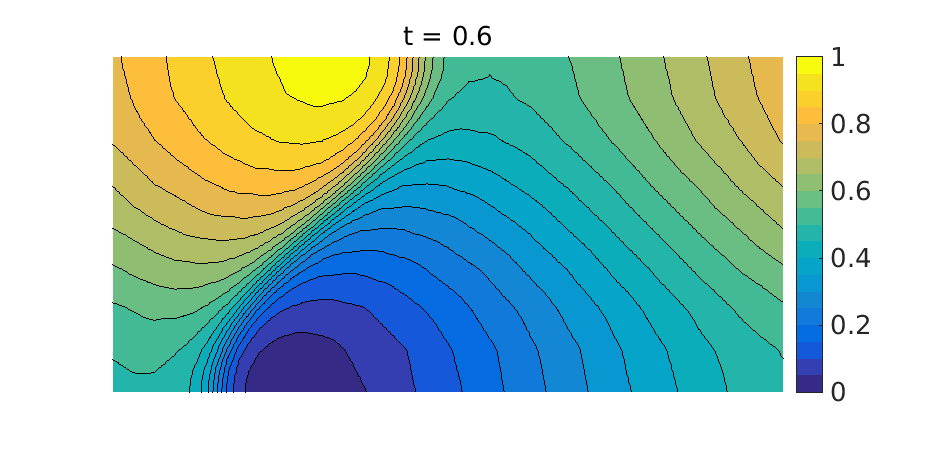}\\  \includegraphics[scale = 0.3,trim=85 25 110 20,clip = true]{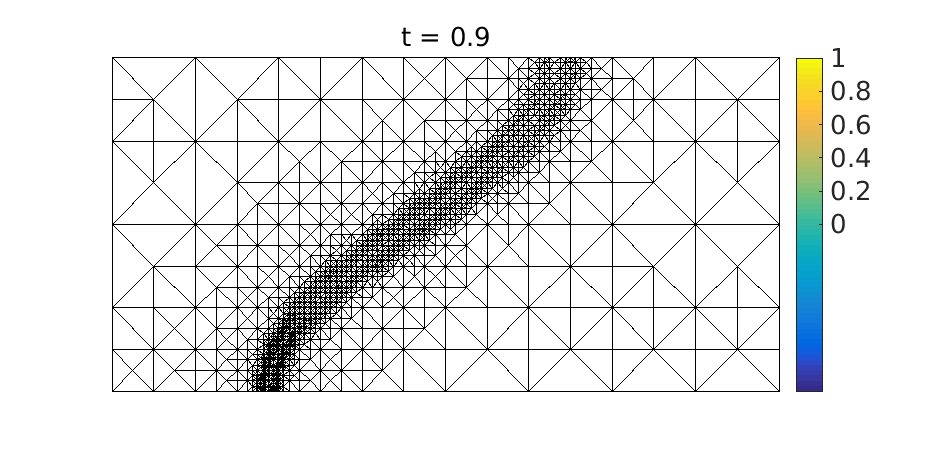}\quad
  \includegraphics[scale = 0.3,trim=85 25  45 20,clip = true]{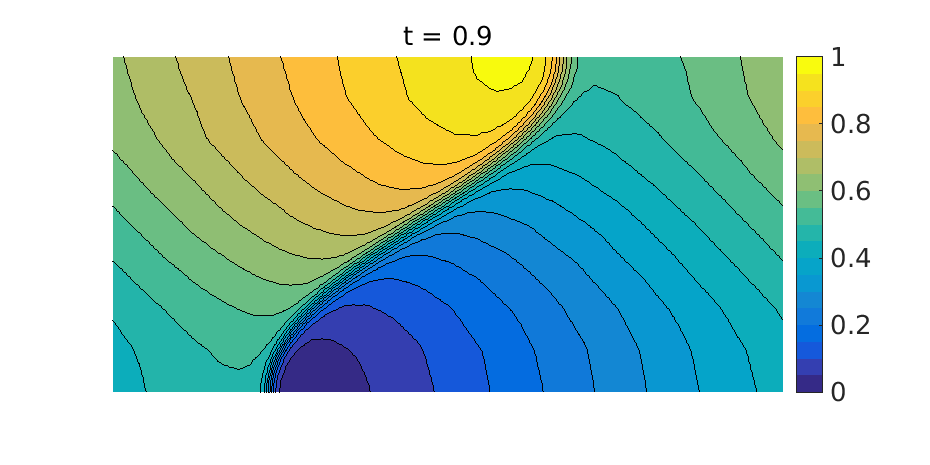}\\  \includegraphics[scale = 0.3,trim=85 25 110 20,clip = true]{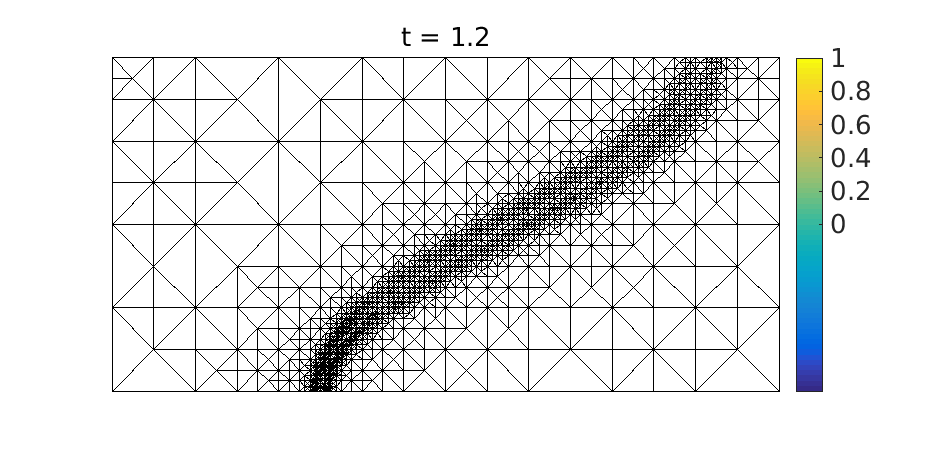}\quad
  \includegraphics[scale = 0.3,trim=85 25  45 20,clip = true]{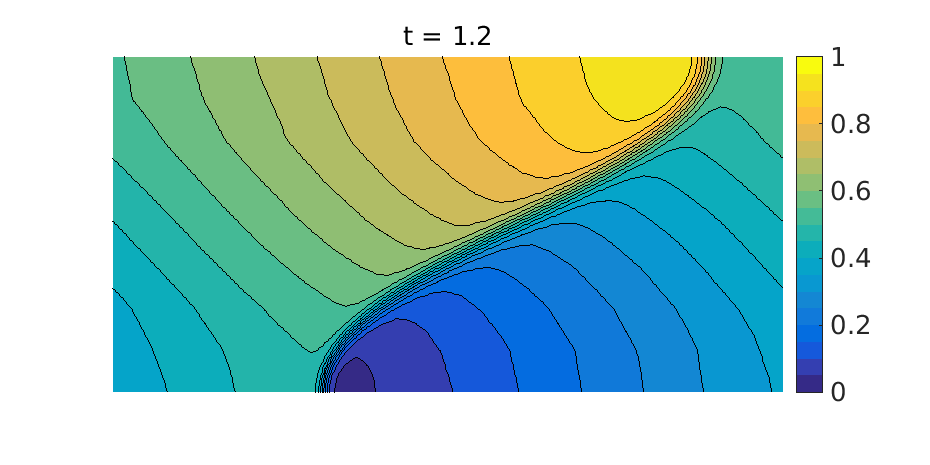}
  \caption{Solution of the Burgers problem with $\parameter=1$ at different times. Left: Adapted finite element meshes, right: numerical solutions}
  \label{fig:nominalSolution}
\end{figure}

\begin{figure}
  \includegraphics[scale = 0.3,trim=85 25 110 20,clip = true]{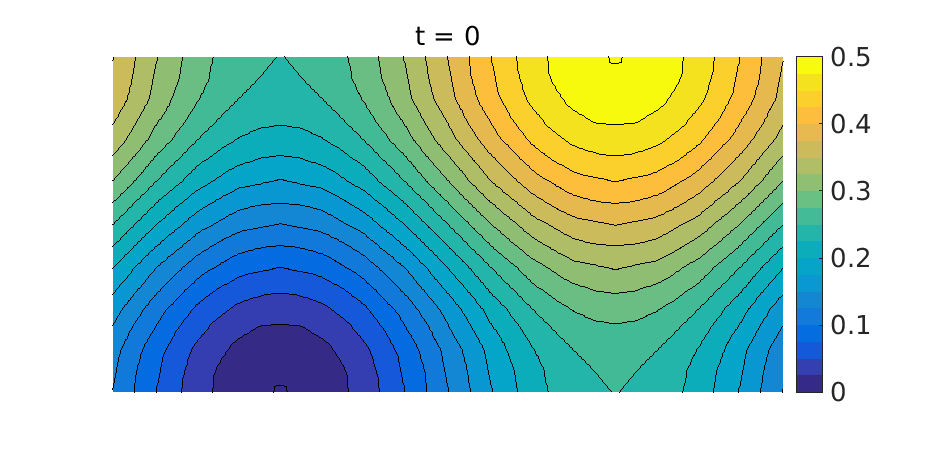}\quad
  \includegraphics[scale = 0.3,trim=85 25  45 20,clip = true]{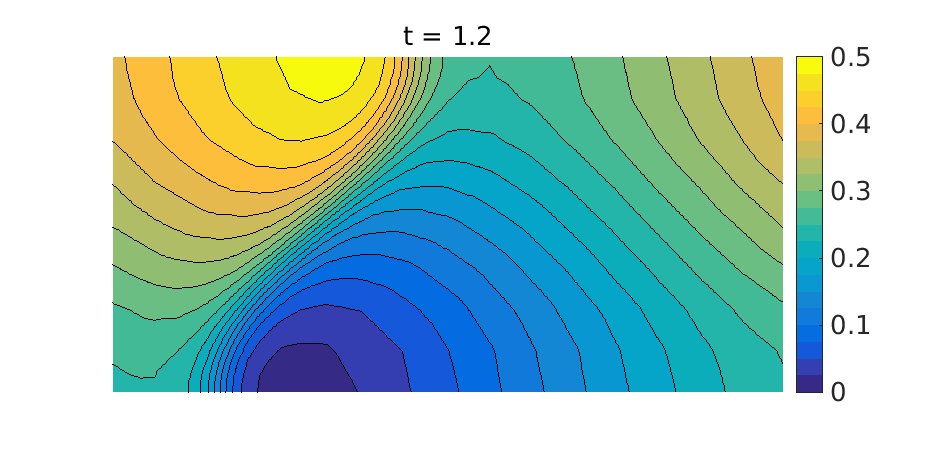}
  \caption{Solution of the Burgers problem with $\parameter=0.5$ at the initial and final times.}
  \label{fig:nominalSolutionB}
\end{figure}

In the presented example, the time discretization uses a step size of $\tau = 0.005$. The spatial discretization in each time step employs the same error estimation, triangle marking and refinement strategies as in \autoref{sec:convectionDiffusionDiscretization}. The adapted meshes and corresponding numerical solutions at various time instances are presented in \autoref{fig:nominalSolution} for a fixed parameter choice of $\parameter = 1$. To illustrate the parameter dependency of the solution, \autoref{fig:nominalSolutionB} shows the respective numerical solutions at the initial and final times for $\parameter = 0.5$.

\subsection{Reduced-order modeling}

{For consistency, we }denote the snapshots by $u_1\in V_1,\dots,u_\nSnap\in V_\nSnap$ with $\nSnap=\nTime\nParam$. We set $u_\iSnap = u_{\iTime,\iParam}$ and $V_\nSnap = V_{\iTime,\iParam}$ with $\nSnap=\iTime+(\iParam-1)\nTime$ for $\iTime=1,\dots,\nTime$ and $\iParam=1,\dots,\nParam$. Using the methods from \autoref{sec:Pod}, we create a POD space $V^\nPod\subset V$ with basis functions $\phi_1,\dots,\phi_\nPod$. Substituting $V$ by $V^\nPod$ in the semi-discretization \eqref{eq:semiTime1}--\eqref{eq:semiTime2} gives rise to the following discretized reduced-order model: Find $u^\nPod_1(\parameter),\dots,u^\nPod_\nTime(\parameter):\parameterDomain\rightarrow V^\nPod$ such that
\begin{alignat}{3}
  (u^\nPod_\iTime-u_0\parameter,\phi)&= 0 &\quad&\forall\phi\in V^\nPod,&\quad&\iTime = 1,\label{eq:rom1}\\
  (u^\nPod_\iTime-u^\nPod_{\iTime-1},\phi) + \tau a(u^\nPod_{\iTime},\phi) + \tau b(u^\nPod_{\iTime},u^\nPod_{\iTime},\phi) &= 0 &\quad&\forall\phi\in V^\nPod,&\quad&\iTime=2,\dots,\nTime.\label{eq:rom2}
\end{alignat}
We can rewrite the model as a set of equations for the POD coefficient vectors: Find $\vec b_1(\parameter),\dots,\vec b_\nTime(\parameter):\parameterDomain\rightarrow \mathbb R^\nPod$ such that
\begin{alignat*}{2}
  M^R\vec b_\iTime - \vec b_0\parameter &= \vec 0, &\quad&\iTime=1,\\
  M^R(\vec b_\iTime-\vec b_{\iTime-1}) + \tau A^R \vec b_{\iTime} + \tau B^R(\vec b_{\iTime})\vec b_{\iTime} &= \vec 0, &\quad&k=2,\dots,\nTime.
\end{alignat*}
Expressions for the constant model coefficient matrices $M^R$, $A^R$, $B^R$ and the initial data $\vec b_0$ follow directly from substitution of the POD expansion \eqref{eq:PODapproximation} into \eqref{eq:rom1}--\eqref{eq:rom2} and testing against the POD basis functions. An evaluation of the term involving $B^R$ amounts to a multiplication with a tensor containing $\nPod^3$ constant model coefficients.

\subsection{Results}

\begin{figure}
  \begin{center}
    \includegraphics[scale = 0.3,trim=10 0 10 0,clip = false]{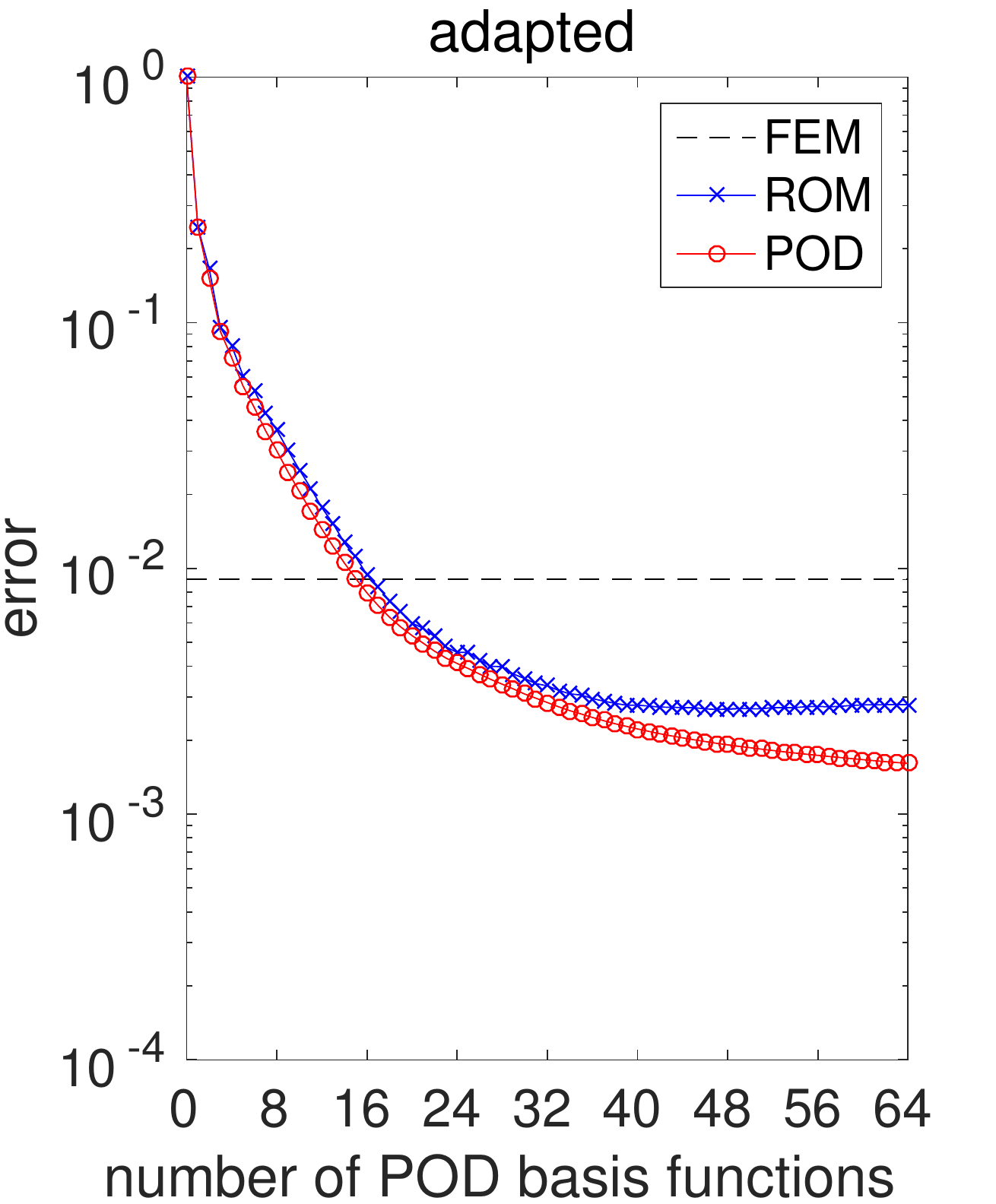}
    \quad
    \includegraphics[scale = 0.3,trim=20 0 10 0,clip = true]{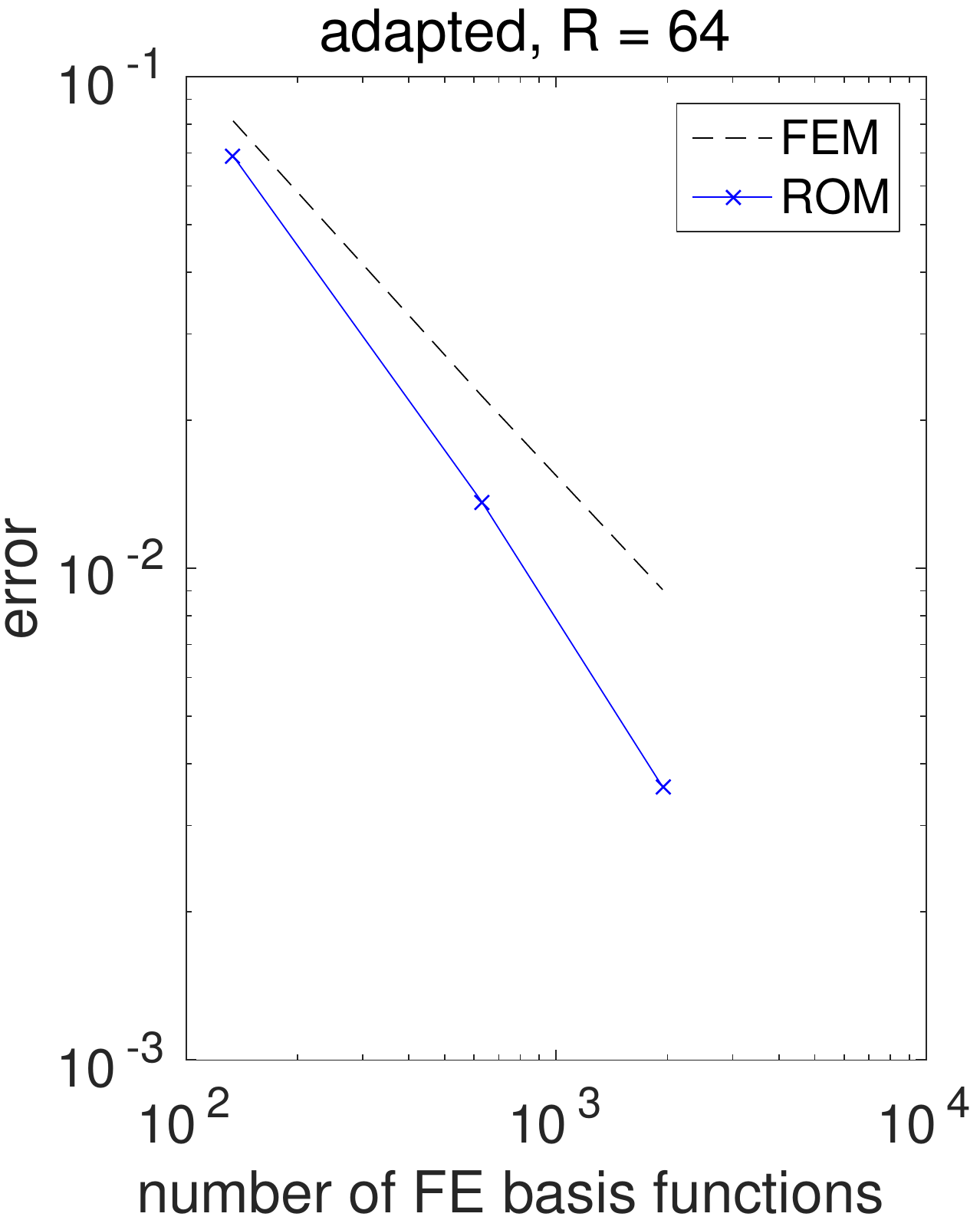}
    \quad
    \includegraphics[scale = 0.3,trim=20 0 10 0,clip = true]{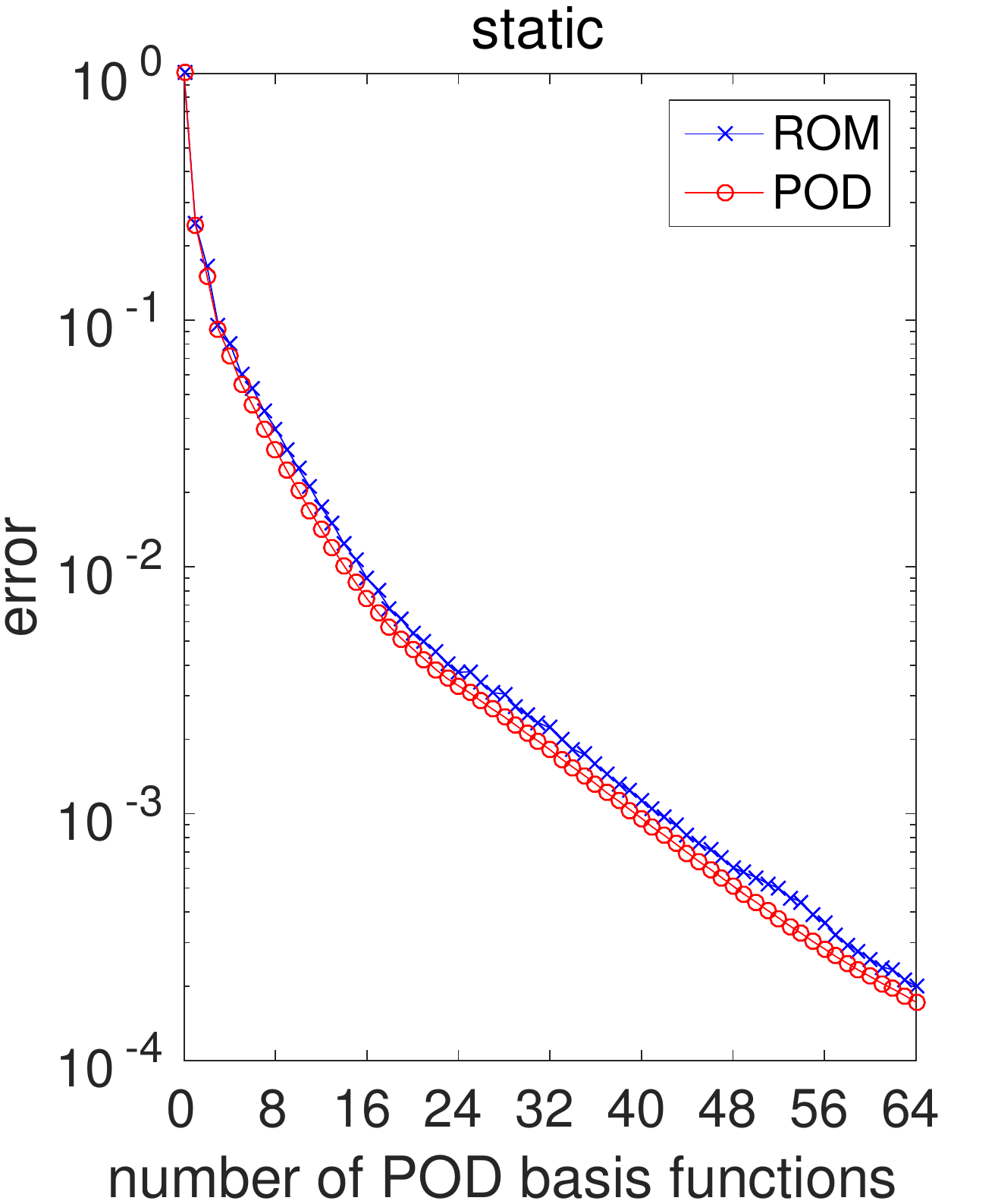}
  \end{center}
  \caption{Error $\epsilon_\text{POD}$ of the POD projection, error $\epsilon_\text{ROM}$ of the solution of the reduced-order model, and finite element discretization error $\epsilon_\text{FEM}$. Left: adapted snapshots, middle: dependence on the finite element refinement for 64 POD basis functions, right: static snapshots computed on the common grid of the space-adapted snapshots.}
  \label{fig:convergence}
\end{figure}

We use the same error measures as in \autoref{sec:convectionDiffusionResults} to test how well the reduced-order models are able to reproduce the snapshots. The results are shown in \autoref{fig:convergence} (left). We observe that the number of POD basis functions needed for a given relative error is increased in comparison with the convection-diffusion problem presented in the previous section. Otherwise the results are similar. In particular, the POD projection error with respect to the snapshots monotonically decreases, but the convergence of the solution of the reduced-order model stagnates at some point. 

\autoref{fig:convergence} (middle) shows the errors depending on the number of finite element degrees of freedom. It can be observed that the level of stagnation is coupled to the finite element discretization error.

The results for a fixed grid, shown in \autoref{fig:convergence} (right), are qualitatively similar to the corresponding results of the convection-diffusion problem, so the same conclusions regarding the use of adaptive snapshots can be drawn in the case of a Burgers problem.

\section{Conclusions}

We have extended the framework of POD-Galerkin reduced-order modeling to snapshot data obtained with adaptive finite elements, where each snapshot may be represented in a different finite element space. The considered reduced-order models rely on a representation in a common finite element space of all snapshots. Because creating such a space may be computationally demanding, we have proposed a method to create the reduced-order model without actually building the common finite element space. 

The POD projection error of our method converges when the POD dimension is increased. The error between the POD Galerkin solution and the snapshots, however, contains a contribution from the spatial finite element discretization, which does not vanish when the POD dimension is increased. We have shown this effect for the case of linear elliptic boundary value problems. 

Our findings are underlined with a numerical example of a convection-diffusion problem. Here we could observe that the error caused by the spatial adaptivity is dominated by the finite element error of the snapshots. Computational results for a Burgers problem suggest that the principal statements can be carried over to a broader class of problems, including non-linear parametrized parabolic PDEs. 

A detailed analysis of the influence of the discretization errors onto the reduced solution for non-linear and/or time-dependent problems with space adapted snapshots is still to be done. The ultimate goal of combining adaptivity with model order reduction is to automate the creation of the reduced-order model by adapting the dimension of the reduced basis and the snapshot discretization according to some global error criterion. To reach this goal, it is still necessary to derive adequate error bounds and to combine them in a global adaptation loop.

\section{Acknowledgments}

This work was supported by the Excellence Initiative of the German Federal and State Governments via the Graduate School of Excellence Computational Engineering at Technische Universit\"at Darmstadt and the Darmstadt Graduate School of Excellence Energy Science and Engineering. The third author was partly supported by the DFG within the collaborative research center TRR154 ``Mathematical Modelling, Simulation and Optimization using the Example of Gas Networks''.


\bibliographystyle{siam}
\bibliography{./adaptive-pod-paper.bib}

\end{document}